\documentclass [a4paper]{article}
\usepackage{amsmath,amsfonts,amsthm}                   %definicje specjalnych fontow np bb (duze tablicowe)
\usepackage{dsfont}
\usepackage{pifont}
\usepackage{graphicx,color}		%grafika
\usepackage{array, tabularx}

\newtheorem{Theorem}{Theorem}
\newtheorem{Corollary}[Theorem]{Corollary}
\newtheorem{Lemma}[Theorem]{Lemma}
\newtheorem{Proposition}[Theorem]{Proposition}

\newtheorem{Definition}{Definition}[section]

\numberwithin{equation}{section}
\numberwithin{Theorem}{section}

\newcommand{\jedynka}{\mathbf{1}{\hskip -3.5 pt}\hbox{l}}

\definecolor{c20}{rgb}{0,0.7,0}
\definecolor{c30}{rgb}{0,0,1}
\definecolor{c40}{rgb}{1,0.1,0.5}
\definecolor{c50}{rgb}{1,0,0}

\begin{document}

\title{Modelling the dependence between a Wiener process and its running maxima and running minima processes}
\author{Piotr Jaworski, Karol D\c{a}browski \\
Institute of Mathematics, University of Warsaw}
\date{\today}

\maketitle

% Abstract (Do not insert blank lines, i.e. \\) 
\abstract{We study a triple of stochastic processes: a Wiener process $W_t$, $t \geq 0$, its running maxima process $M_t=\sup \{W_s: s \in [0,t]\}$, and its running minima process $m_t=\inf \{W_s: s \in [0,t]\}$. We derive the analytical formula for the  corresponding copula and show that it is supported on the hemicube, a convex hexahedron with seven vertices.
As an application, we draw out an analytical formula for pricing of a double barrier option.}

%
%\abstract{We study a triple of stochastic processes: a Wiener process $W_t$, $t \geq 0$, its running maxima process $M_t=\sup \{W_s: s \in [0,t]\}$ and its running minima process $m_t=\inf \{W_s: s \in [0,t]\}$. We derive the analytical formula for the  corresponding copula and show that it is supported on the hemicube, a convex hexahedron with 7 vertices.
%As an application we draw out an analytical formula for pricing of a double barrier option.}
%
%\noindent
%Keywords: Copulas; Wiener process; Running maxima and minima; Strong Markov property; Reflection principle;  Spearman rho; Double-barrier options.\\
%MSC2020: 62H05, 60G18, 60G70, 60E05, 62H20, 91G20.\\

\section{Introduction}\label{S:Introduction}

The Wiener process is one of the most ubiquitous stochastic processes
used to model complex phenomena. One of such areas is
financial mathematics, where this process is the foundation for all
models of stock, derivatives, or portfolio analysis. But with
the development of the market, more complex derivatives have shown up,
and a simple approach to the Wiener process was no longer enough. In particular, there was a
question as to what is the distribution of the Wiener process adding
boundaries on trajectories. It is especially important when one is dealing with
barrier options.

The structure of article is as follows.
In the first part, there will be provided the main results. We derive the
 cumulative distribution function and copula for a triple consisting of a Wiener process
and its supremum and infimum. 
{{The} 
 cumulative distribution function is supported on a cone and the copula on a hemicube.}
Then, we present special cases of this distribution, particularly focusing on marginals.
The next part is devoted to proofs of stated theorems. 
The starting point is the
 strong Markov property of Wiener processes. First, we study the simple case based on one stopping time. Finally, to the obtain desired results, we move to double stopping time. We provide a recursive dependency, which will lead to the characterization of the cumulative distribution the function and copula in terms of series.

The last parts show applications in financial mathematics, i.e., we calculate the price of a  European double-barrier option in the Black--Scholes model.
{{Compare}~\cite{Cherubini_Romagnoli_2010, Haug_1999, Barker_2007}.  Our approach is based on the Girsanov theorem, which implies the existence a probabilistic measure $Q$ under which the increment of the logarithm of a stock price process divided by the volatility follows the Wiener law. The double change of measure, from the risk neutral to $Q$ and back after some integral calculations, significantly simplifies the proof.} 

\section{Main Results}\label{S:Main}

Consider a probability space $(\Omega, \mathcal{F}, \mathbb{P})$ and a Wiener process $(W_t)_{t\ge0}$ with natural filtration with respect to that process, i.e., $\mathcal{F}_t = \sigma \left\{\left.W_{{s}}^{{-1}}(A)\right|s \le t,,A\in \mathcal{B}(\mathbb{R}) \right\}$---where $\mathcal{B}(\mathbb{R})$ is the Borel $\sigma$-algebra on $\mathbb{R}$.

{We denote as} %%%: 1. Please recheck all equations and make sure there are no duplicated equations in the whole manuscript. Thanks! 
%2. Please carefully check variable formatting (italic, bold, subscript, uppercase, etc.) throughout the manuscript to ensure the formatting is consistent and revise if needed..
$$M_t := \sup_{0 \le s \le t}W_s,$$
and 
$$m_t := \inf_{0 \le s \le t}W_s$$
the associated running maxima and minima processes.

Note that the triple $(W_t,M_t,m_t)$, $t> 0$ is a self-similar {process:} %%%: Equations should be numbered in order without section number, we revised all, please confirm
\begin{equation}
\forall t>0 \;\;\; (W_t,M_t,m_t) \stackrel{d}{=} (\sqrt{t}W_1, \sqrt{t}M_1, \sqrt{t}m_t).
\end{equation}

{Let $F_t(x,y,z)$ be the cumulative distribution function of the joint distribution $(W_t,M_t,m_t)$}, where $t>0$.
As follows from the celebrated Sklar theorem,
there exists a copula process, $C_t$, $t>0$, such that
\[F_t(x,y,z) = C_t\big( \mathbb{P}(W_t\leq x), \mathbb{P}(M_t\leq y), \mathbb{P}(m_t\leq z) \big).\]

{For} %%%: The sentence or paragraph which start as capitalized letter below the single-line formula should be indented. So we added the indent. Please confirm all
 more details about copula theory and some of its applications, 
we refer to
~\cite{Nel, Durante_Sempi_2016, 10CTA, 13CMQ, Joe97, Joe15, CLV}.
Since the process $(W_t,M_ t,m_t)$ is self-similar,
 the copula $C_t$ of the triple $(W_t,M_t,m_t)$, where $t> 0$, does not depend on $t$. We will denote it as $C_{W,M,m}$.
 
Based on the strong Markov property of a Wiener process and the Sklar theorem, we prove what follows, where $\Phi$ denotes the cumulative distribution
of the $N(0,1)$ probability~law.
\begin{Theorem}\label{T:C}
{The joint} %%%: Definitions/theorems/remarks/lemma, etc., should be numbered sequentially. We revised all accordingly, please confirm
 cumulative distribution function $F_t(x,y,z)$ of $(W_t, M_t,m_t)$, where $t>0$, is of the~form
$$ F_t(x,y,z) =C_{W,M,m}\left( \Phi(\frac{x}{\sqrt{t}}), \left( 2\Phi(\frac{y}{\sqrt{t}}) -1\right)^+, \min(2\Phi(\frac{z}{\sqrt{t}}),1) \right), \;\;\; t >0,$$
where
%%\begin{adjustwidth}{-\extralength}{0cm}
%\centering %% If there is a figure in wide page, please release command \centering
\[C_{W,M,m}(u,v,w)= \left\{
\begin{array}{cl}
 u-\Phi\big(\Phi^{-1}(u)-r\big) & \mbox{for }   2u\leq w,  \\
 2\Psi\big(\Phi^{-1}(w/2),r,s\big) - \Psi\big(\Phi^{-1}(u)-s,r,s\big) &\\
-\Psi\big(-\Phi^{-1}(u)-s+r,r,s\big) & \mbox{for }   w < 2u \leq 1+v,\\
2\Psi\big(\Phi^{-1}(w/2),r,s)\big) - 2\Psi\big(2\Phi^{-1}(w/2)-\Phi^{-1}((1+v)/2),r,s\big)  &\mbox{for }   1+v < 2u ,
\end{array}
\right.\]
%%%\end{adjustwidth}
where $\Psi$ is given by formula 
\begin{equation}\label{E:Psi}
\Psi(q,r,s) = \sum_{k=0}^\infty (\Phi(q-ks) - \Phi(q-r-ks)) 
\end{equation}
with
\[ r=2\Phi^{-1}\big((1+v)/2\big) \;\; \mbox{and } \;\; s=2\big(\Phi^{-1}((1+v)/2)-\Phi^{-1}(w/2))\big).\]
\end{Theorem}

\hspace*{5mm}

The proof is provided in Section \ref{S:Proof_2}.

Note that the copula $C_{W,M,m}$ is absolutely continuous with respect to the Lebesgue measure and is supported on the polyhedron
\[ P=\{(u,v,w) \in [0,1]^3:\;  w \leq 2u \leq 1+v \},\]
which has seven vertices
\[(0,0,0),\; (0,1,0),\; (1,1,0),\; (1,1,1),\; (1/2,0,1), \; (1/2,0,0), \; (1/2,1,1),\]
six faces, two rectangles, two trapezoids, two triangles, and eleven edges. 
Such polyhedrons are called hemicubes because they are linearly equivalent to a polyhedron obtained by cutting the unit cube by half
by the surface containing two opposite vertices (for example, $(0,0,0)$ and $(1,1,1)$) ({Figure} %%%: All Figure should be cited, we add the citation of Figure 1 here; if it is not appropriate, please modify and ensure that the first citation of each Figure appears in numerical order..
 \ref{R:P}).%\vspace{-40pt}

\begin{figure}[ht]
%\begin{center}
\begin{picture}(200,200) (0,0)

{
\put(10,10){\circle*{3}}
\put(10,110){\circle*{3}}
\put(110,110){\circle*{3}}
\put(150,160){\circle*{3}}
\put(100,60){\circle*{3}}
\put(100,160){\circle*{3}}
\put(60,10){\circle*{3}}
}
\put(110,10){\circle*{3}}
\put(50,60){\circle*{3}}
%\put(50,160){\circle{3}}
%\put(150,60){\circle{3}}

{
\put(10,10){\line(0,1){100}}
\put(10,110){\line(1,0){100}}
%\put(110,110){\line(1,1){50}}
\put(10,10){\line(1,0){50}}
%\put(10,10){\line(2,1){100}}
%\put(60,10){\line(1,1){50}}
\put(60,10){\line(1,2){50}}
\put(100,60){\line(1,2){50}}
\put(100,60){\line(0,1){100}}
%\put(10,110){\line(2,1){100}}
\put(100,160){\line(1,0){50}}
\qbezier(10,10)(55,35)(100,60)
\qbezier(60,10)(80,35)(100,60)
\qbezier(10,110)(55,135)(100,160)
\qbezier(110,110)(130,135)(150,160)
}

\put(10,110){\vector(0,1){30}}
\put(60,10){\vector(1,0){80}}
\put(10,10){\line(4,5){50}}
\put(60,72){\vector(3,4){1}}
\put(143,10){u}
\put(10,140){v}
\put(60,72){w}
\put(10,0){0}
\put(56,0){1/2}
\put(108,0){1}
\put(2,105){1}
\put(42,56){1}
\end{picture}

%\end{center}
\caption{{Polyhedron} %%%: 1. we revised the color of line into black in image, please confirm 2
 $P${---}the support of the $C_{W,M,m}$ copula.}\label{R:P}
\end{figure}

Furthermore, the copula $C_{W,M,m}$ is invariant with respect to the involution induced by the rotation of the unit cube by an angle $\pi$ about the line 
\[ \{(u,v,w): u=0.5, v+w=1 \}. \]

For every triple of generators $(U,V,W)$ of the copula $C_{W,M,m}$, the triple $(1-U,1-W,1-V)$ is a set of generators of this copula as well. Hence,%\vspace{-12pt}
%%\begin{adjustwidth}{-\extralength}{0cm}
%\centering %% If there is a figure in wide page, please release command \centering
\begin{eqnarray}\label{E:Involution_1}
C_{W,M,m}(u,v,w)&=&\mathbb{P}(U\leq u, V \leq v, W \leq w)= \mathbb{P}(1-U\leq u, 1-W \leq v, 1-V \leq w) \\
&=& u+v+w -2 + C_{W,M,m}(1,1-w,1-v) +C_{W,M,m}(1-u,1,1-v)  \nonumber \\
 &&+ \, C_{W,M,m}(1-u,1-w,1) - C_{W,M,m}(1-u,1-w,1-v) .\nonumber
\end{eqnarray}
%%%\end{adjustwidth}

The marginal copulas are given by %\vspace{-12pt}
%%\begin{adjustwidth}{-\extralength}{0cm}
%\centering %% If there is a figure in wide page, please release command \centering
\begin{eqnarray*}
C_{W,M}(u,v)=C_{W,M,m}(u,v,1)&=& \left\{
\begin{array}{cl}
u-\Phi\big(\Phi^{-1}(u) -2\Phi^{-1}((1+v)/2)\big)  & \mbox{for } 2u \leq 1+v,\\
v & \mbox{for }  1+v < 2u,
\end{array} \right.\\
C_{W,m}(u,w)=C_{W,M,m}(u,1,w)&=& \left\{
\begin{array}{cl}
u  & \mbox{for } 2u \leq w,\\
w- \Phi\big(2\Phi^{-1}(w/2) - \Phi^{-1}(u)\big) \;\;\;\; & \mbox{for } w <2u,
\end{array} \right.\\
C_{M,m}(v,w)=C_{W,M,m}(1,v,w)&=& 
2 \Psi\big(\Phi^{-1}(w/2),r,s) - 2 \Psi\big(2\Phi^{-1}(w/2)-\Phi^{-1}((1+v)/2),r,s\big).
\end{eqnarray*}
%%%\end{adjustwidth}

All three of the above bivariate copulas are absolutely continuous. 
The first two are supported on trapezoids.
The first one is supported on %Please check meaning retained.
\[ \{(u,v) \in [0,1]^2 : 2u \leq 1+v \},\]
the second one is supported on 
\[ \{(u,w)\in [0,1]^2 : w \leq 2u \}.\]

\begin{figure}[ht]
%\begin{center}
 \includegraphics[width=13.5cm]{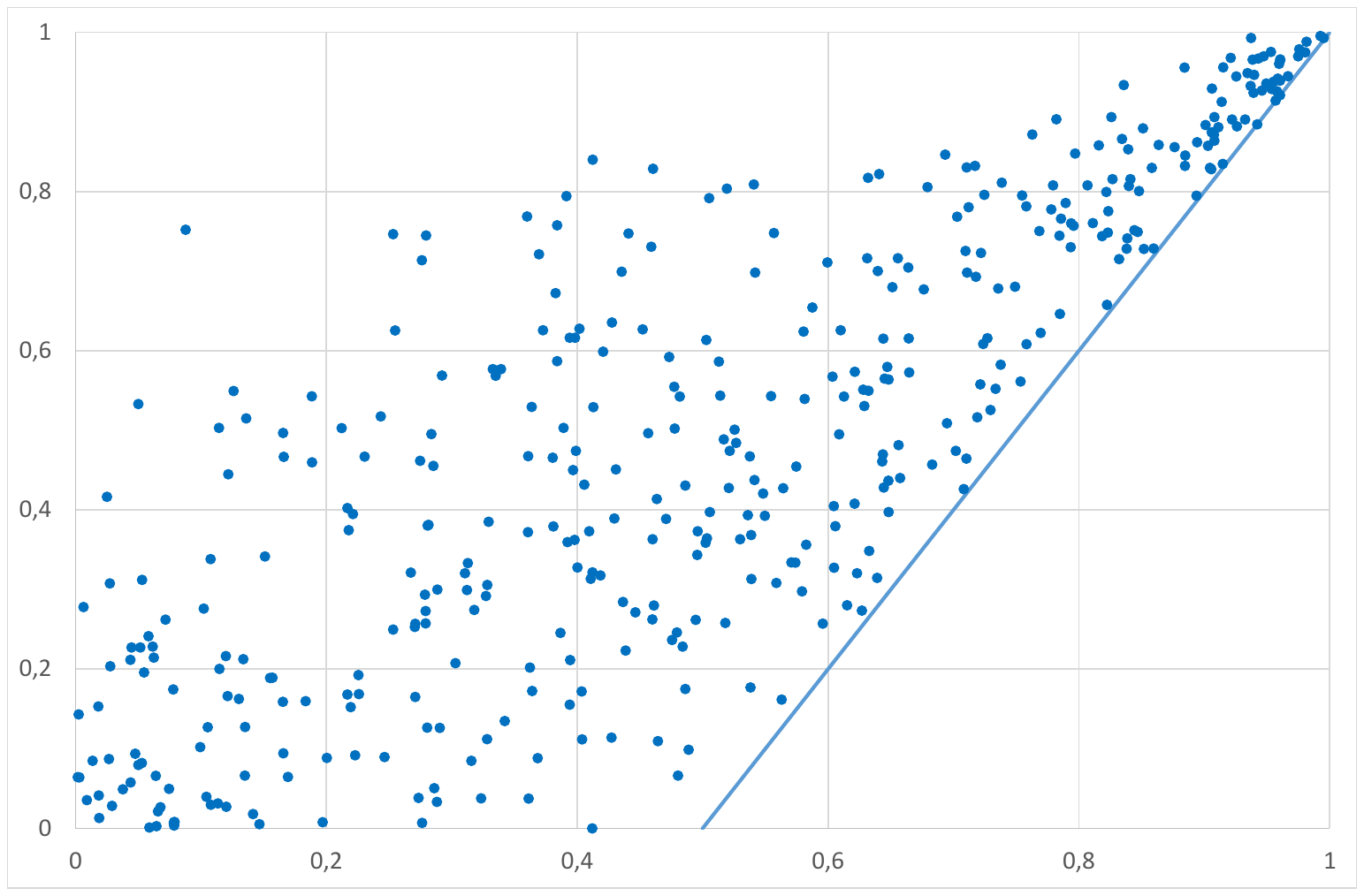}
% \end{center}
\caption{{Copula} %%%: Please use periods as decimal signs instead of commas, e.g., "0,1" should be "0.1".
 $C_{W,M}$ (scatterplot).}\label{R:WM}
\end{figure}

\begin{figure}[ht]
%\begin{center}
 \includegraphics[width=13.5cm]{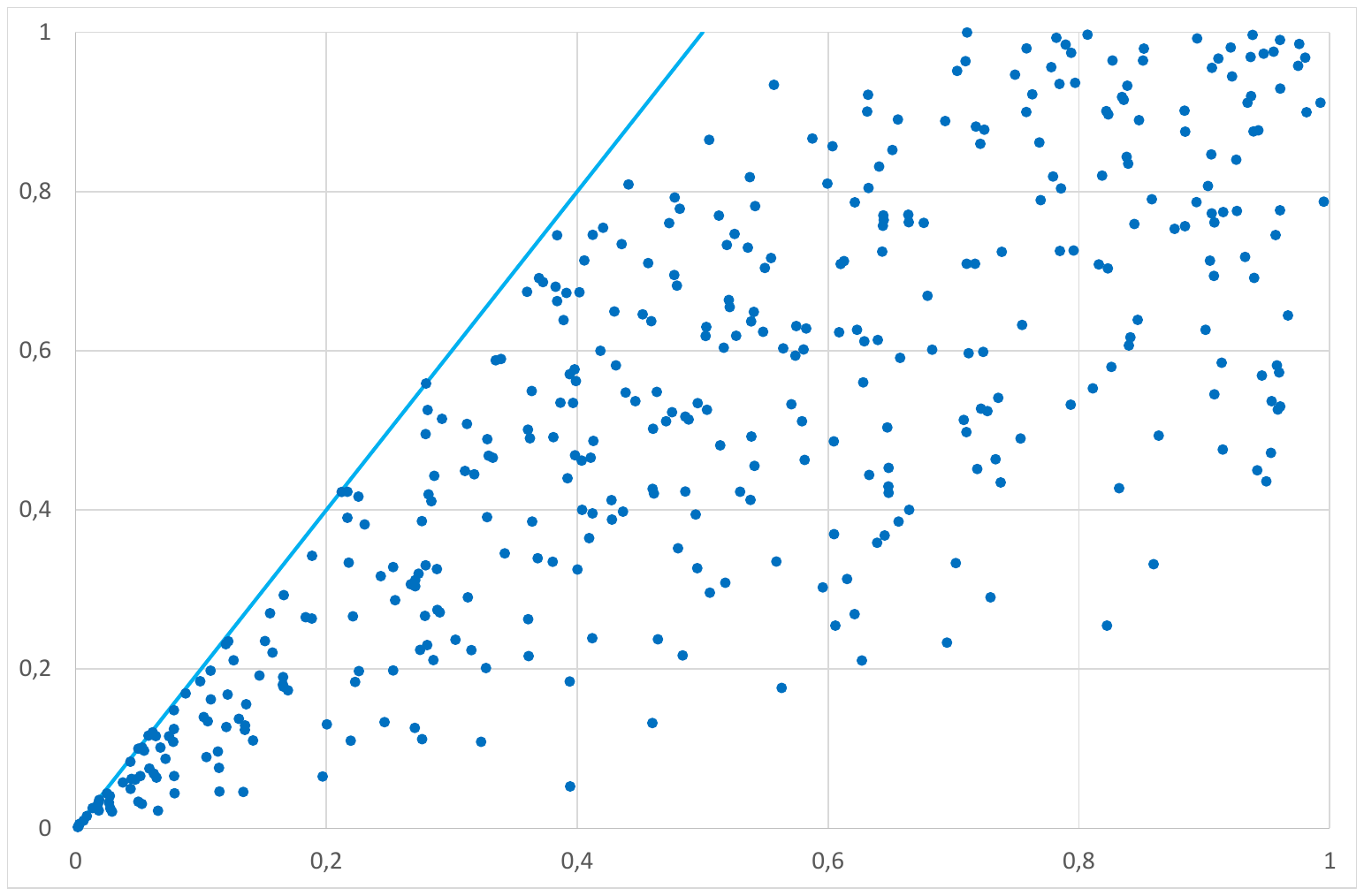}
% \end{center}
\caption{{Copula} %%%: Please use periods as decimal signs instead of commas, e.g., "0,1" should be "0.1".
 $C_{W,m}$ (scatterplot).}\label{R:Wm}
\end{figure}

\begin{figure}[ht]
%\begin{center}
 \includegraphics[width=13.5cm]{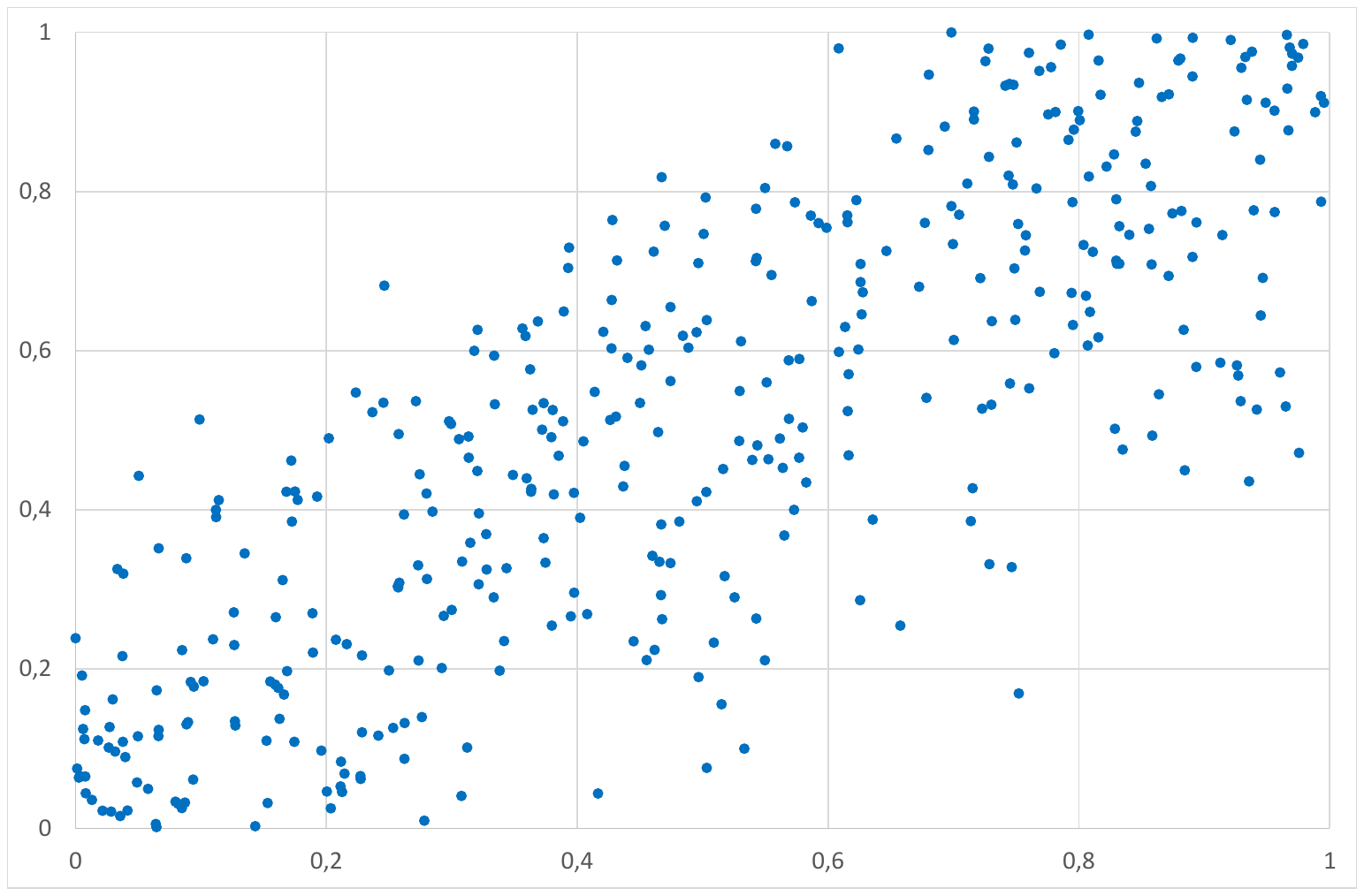}
% \end{center}
\caption{{Copula} %%%: 
 $C_{M,m}$ (scatterplot).}\label{R:Mm}
\end{figure}

The third one has a full support $[0,1]^2$. %Compare 
{The scatterplots of the above copulas are shown on Figures \ref{R:WM}--\ref{R:Mm}.}

 %\vspace{-6pt}

Furthermore, the first and second copula are  duals of each other: $C_{W,m}$ is  a survival copula of $C_{W,M}$ (and vice versa). 
As a matter of fact, they are reflections   of a copula introduced in~\cite{Jaworski_2023}. See also~\cite{Ades_atal_2022}.
We have
\begin{equation}\label{E:Involution_2}
 C_{W,m}(u,w) = u+w-1 +C_{W,M}(1-u,1-w).
 \end{equation}
 
The third one is self-dual with respect to reflection about the line $v+w=1$,
\begin{equation}\label{E:Involution_3}
 C_{M,m}(v,w) = u+v-1 +C_{M,m}(1-w,1-v).
 \end{equation}\vspace{6pt}

All three marginal copulas are highly dependent.

\begin{Proposition}\label{P:rho}
{The} 
 Spearman $\rho$  of the marginal copulas of the copula $C_{W,M,m}$ is given {by} %%%: 

\begin{eqnarray*}
 \rho(C_{W,M})&=&\rho(C_{W,m} ) \;\; = 2 - \frac{6}{\pi} \arccos\left( \frac{\sqrt{6}}{3} \right) \approx 0.8245,\\
\rho(C_{M,m})&=& \frac{6}{\pi} \sum_{n=0}^\infty (-1)^n \bigg(
\arccos \left( \frac{n}{\sqrt{2(n(n+1)+1)}}\right)
+ \arccos \left( \frac{n+1}{\sqrt{2(n(n+1)+1)}}\right) \nonumber \\
&& - \, \arccos \left( \frac{n+1}{\sqrt{2((n+1)(n+2)+1)}}\right)
- \arccos \left( \frac{n+2}{\sqrt{2((n+1)(n+2)+1)}}\right)
\bigg) \\
&\approx & 0.80649. \nonumber
\end{eqnarray*}
%%\end{adjustwidth}
\end{Proposition}

The proofs are provided in Section \ref{S:Properties}.

\section{Proofs and Auxiliary Results}\label{S:Proofs}
\subsection{Strong Markov Property}\label{S:Markov}

The strong Markov property plays a crucial role in the proof of Theorem \ref{T:C}.
We recall the definition.
\begin{Definition}
{Suppose} 
 that  $X=(X_{t}:t\geq 0) $ is a stochastic process on a probability space  $ (\Omega ,{\mathcal {F}},\mathbb{P}) $ 
with natural filtration  $ \{{\mathcal {F}}_{t}\}_{{t\geq 0}}$. Then, for any stopping time $ \tau$   on  $ \Omega$ , we can define 
\[  {{\mathcal {F}}_{\tau }=\{A\in {\mathcal {F}}:\forall t\geq 0,\{\tau \leq t\}\cap A\in {\mathcal {F}}_{t}\}}.\]

Then,  $ X$ is said to have the strong Markov property if, for each finite stopping time  $ \tau$ ,
 we have that for each  $ t\geq 0$,  $X_{{\tau +t}} - X_\tau$ is independent of  ${\mathcal {F}}_{\tau }$.
\end{Definition} 

It is well known that the Wiener process fulfills the above definition; see, for example,~\cite{Billingsley_1986} ({Section 37,} 
 Theorem 37.5), \cite{Baudoin_2014} (Section 3.2), or~\cite{Schilling_2014} (Theorem 6.5). Furthermore, for a Wiener process $W_t, t \geq 0$, we have 

\begin{Theorem}\label{T:SMP}
Let $\tau$ be a finite stopping time and put
\begin{equation}
W^\ast_t(\omega)= W_{\tau+t}(\omega) - W_\tau(\omega), \;\;\; t \geq 0, \;\; \omega \in \Omega.
\end{equation}

Then, $W^*_t: t \geq 0$ is a Brownian motion, and it is independent of  ${\mathcal {F}}_{\tau }$.
\end{Theorem}

Based on the fact that $W^\ast_t$ is symmetric
\[ \forall t>0 \;\;\; W^\ast_t\stackrel{d}{=} -W^\ast_t \]
we obtain the following version of the reflection principle.

\begin{Theorem}
Let $\tau$ be a finite stopping time and let a set $A$ be an element of ${\mathcal {F}}_{\tau }$ of positive probability on which  $\tau \leq t$, 
\[ A \in {\mathcal {F}}_{\tau }, \;\; \mathbb{P}(A)>0, \;\; A \subset \{\omega : \tau(\omega) \leq t \},\] 
{such that} under the condition  $A$, the probability laws of $W_t$ and $2 W_\tau - W_t$ coincide:
\begin{equation}
W_t  \, | \,  A \; \stackrel{d}{=} 2 W_\tau - W_t  \, | \,   A .
\end{equation}
\end{Theorem}

%\noindent
\begin{proof}
{We} 
 observe that under the condition $\tau \leq t$,  
\begin{eqnarray}
W_t &=& W_t - W_{\tau} + W_\tau = W^\ast_{t-\tau} +W_\tau= W^\ast_{|t-\tau|} +W_\tau ,  \\
2 W_\tau - W_t &=& - (W_t -  W_{\tau}) + W_\tau= - W^\ast_{t-\tau} +W_\tau = - W^\ast_{|t-\tau|} +W_\tau.
\end{eqnarray}

Next, we compare the conditional distribution functions of the right sides of the above equations.
We apply the fact that  the random variables $\tau$ and $ W_\tau$ are  $\mathcal{F}_\tau$ measurable {and} 
\begin{eqnarray}
 \mathbb{P}( W^\ast_{|t-\tau|} \leq  x  - W_\tau \, | \, \mathcal{F}_\tau) 
&\stackrel{{a.s.}}{=}&
\left\{\begin{array}{cl}
\jedynka_{W_\tau \leq x} & \mbox{if } \tau=t\\
  \Phi\left( \frac{x-W_\tau}{\sqrt{t-\tau}} \right)  & \mbox{other-ways }
\end{array} \right. \\
&=&
\left\{\begin{array}{cl}
\jedynka_{W_\tau \leq x} & \mbox{if } \tau=t\\
  1 -\Phi\left( \frac{W_\tau -x}{\sqrt{|t-\tau|} }\right)  & \mbox{other ways }
\end{array} \right. \\
&\stackrel{{a.s.}}{=}&  \mathbb{P}( -W^\ast_{|t-\tau|} \leq  x  - W_\tau \, | \, \mathcal{F}_\tau) \nonumber 
 \end{eqnarray}
where $\Phi$ denotes the cumulative distribution function of a standard normal random variable ($N(0,1)$).

To conclude the proof, we apply the second point of the definition of the conditional expected value (see~\cite{Billingsley_1986} \S 34). %\vspace{-6pt}
{\small
\begin{eqnarray}
\mathbb{P}( W^\ast_{t-\tau} +W_\tau  \leq x \, | \, A) &=&
\frac{\mathbb{E} \left( \jedynka_{W^\ast_{|t-\tau|} +W_\tau  \leq x} \jedynka_A \right)}{\mathbb{P}(A)}\\
&=& \frac{1}{\mathbb{P}(A)} \int_A \jedynka_{W^\ast_{|t-\tau|} +W_\tau  \leq x} d\mathbb{P} \nonumber \\
&=& \frac{1}{\mathbb{P}(A)} \int_A \mathbb{E}\left(\jedynka_{W^\ast_{|t-\tau|} +W_\tau  \leq x} \; |\;  \mathcal{F}_\tau \right) d\mathbb{P} \nonumber \\
&=& \frac{1}{\mathbb{P}(A)} \int_A \mathbb{E}\left(\jedynka_{-W^\ast_{|t-\tau|} \leq  x  - W_\tau} \; |\;  \mathcal{F}_\tau \right) d\mathbb{P} \nonumber \\
&=& \frac{1}{\mathbb{P}(A)} \int_A \jedynka_{-W^\ast_{|t-\tau|} \leq  x  - W_\tau}  d\mathbb{P} \nonumber \\
&=& \frac{\mathbb{E} \left( \jedynka_{- W^\ast_{|t-\tau|} +W_\tau  \leq x} \jedynka_A \right)}{\mathbb{P}(A)} \nonumber\\
&=& \mathbb{P}( -W^\ast_{t-\tau} +W_\tau  \leq x \, | \, A) . \nonumber
\end{eqnarray}}
\end{proof}

\subsection{Marginal Distributions}\label{S:Margins}

We consider two families of stopping times. One corresponds to an upper  threshold, and the second corresponds to a lower one. In more details,
for $y>0$ and $z<0$, we put
\begin{eqnarray}
T_{y} &=& \inf \{t> 0: W_t = y\},\\
T_{z} &=& \inf \{t> 0: W_t = z\}.
\end{eqnarray}

Note that for fixed positive $t$, the event $T_y \leq t$ is equivalent to the event $M_t \geq y$:
\begin{equation}
\{\omega : T_y(\omega) \leq t\} = \{\omega : M_t(\omega) \geq y \}
\end{equation}
and the event  $T_z \leq t$ is equivalent to the event $m_t \leq z$:
\begin{equation}
\{\omega : T_z(\omega)  \leq t\} = \{\omega : m_t(\omega) \leq z \}.
\end{equation}

We apply the reflection principle for $\tau=T_y$, where $y>0$. 
Since 
\[W_{T_y}=y,\]
the reflection principle simplifies to
\begin{equation}
W_t \; | T_y  \leq t \stackrel{d}{=} 2y -W_t | T_y  \leq t .
\end{equation}

For  $y,t >0$ we obtain
%\vspace{-3pt}
\begin{eqnarray}
\mathbb{P}(W_t \leq x, M_t \geq y)&=&
\mathbb{P}(W_t \leq x, T_y \leq t)\\
&=& \mathbb{P}(W_t \leq x\; | T_y \leq t) \mathbb{P}( T_y \leq t) \nonumber \\
&=& \mathbb{P}(2y-W_t \leq x\; | T_y \leq t) \mathbb{P}( T_y \leq t) \nonumber \\
&=& \mathbb{P}(W_t \geq 2y- x, M_t \geq y). \nonumber
\end{eqnarray}

The further assumption $y \geq x$ implies that the event $M_t \geq y$ contains the event $W_t \geq 2y- x$ 
\[ \{\omega: W_t(\omega) \geq 2y- x \} \subset \{\omega : M_t(\omega) \geq y\}.\]

Indeed, for every $\omega $ from the first set
\[M_t (\omega) \geq W_t(\omega) \geq 2y-x \geq y.\]

Therefore, for $y \geq x$ and $y,t>0$,
\begin{eqnarray}\label{E:WM}
\mathbb{P}(W_t \leq x, M_t \geq y)&=&
\mathbb{P}(W_t \geq 2y- x, M_t \geq y)\\
&=& \mathbb{P}(W_t \geq 2y- x) = 1- \Phi\left(\frac{2y-x}{\sqrt{t}}\right)= \Phi\left(\frac{x-2y}{\sqrt{t}}\right), \nonumber
\end{eqnarray}
where $\Phi$ denotes the cumulative distribution function of a standard normal random variable ($N(0,1)$).

Similarly, for $t>0$, $z<0$, and $z \leq x$, putting $\tau=T_z$,  we obtain
\begin{eqnarray}\label{E:Wm}
\mathbb{P}(W_t \geq x, m_t \leq z)&=&
\mathbb{P}(W_t \leq 2z- x, m_t \leq z)\\
&=& \mathbb{P}(W_t \leq 2z- x) = \Phi\left(\frac{2z-x}{\sqrt{t}}\right). \nonumber
\end{eqnarray}

The above, together with inequalities
\[ m_t \leq W_t \leq M_t, \;\; m_t \leq 0 \leq M_t,\]
allows us to determine the joint distribution functions for the pairs $(W_t,M_t)$ and $(W_t,m_t)$.

\begin{Theorem}\label{T:F_WM}
The joint cumulative distribution function $F_{W,M,t}(x,y)$ of $(W_t,M_t)$, where $t>0$, is of the form
$$ F_{W,M,t}(x,y) =F_{W,M}\left( \frac{x}{\sqrt{t}},\frac{y}{\sqrt{t}}\right), \;\;\; t >0,$$
where
$$F_{W,M}(x,y)= \left\{
\begin{array}{ccc}
0 & \mbox{for } & y < 0,\\
2 \Phi(y) -1 & \mbox{for } & 0 \leq y < x,\\
\Phi(x) - \Phi(x-2y) & \mbox{for } &  x^+ \leq y .
\end{array}
\right. $$

Furthermore,
\[ \mathbb{P} (M_t \leq y) = \max\left( 0, 2\Phi\left( \frac{y}{\sqrt{t}} \right) -1 \right).\] 
\end{Theorem}

%\noindent
\begin{proof}
Case $y<0$.

Since $M_t \geq 0$, 
\[ \forall t >0\;\;\; \mathbb{P}(M_t <0)=0.\]

%\noindent
Case $y\geq \max(0,x)$.

We apply Formula (\ref{E:WM}).
\begin{eqnarray}
\mathbb{P}(W_t\leq x, M_t < y) &=& \mathbb{P}(W_t\leq x) - \mathbb{P}(W_t\leq x, M_t \geq y)\\
&=& \Phi\left( \frac{x}{\sqrt{t}} \right) - \Phi\left( \frac{x-2y}{\sqrt{t}}\right). \nonumber
\end{eqnarray}

Due to the continuity of $\Phi$, we obtain
\[F_{WM}(x,y) =\Phi(x) - \Phi(x-2y).\]
  
%\noindent
Case $x=y$.

The event $\{\omega : M_t(\omega) \leq y\}$ is contained in the event $\{\omega : W_t(\omega)\leq y\}$. Indeed, for $\omega$ belonging to the first set, we have
\[ W_t(\omega) \leq M_t(\omega) \leq y .\]

Hence, from the previous step, we obtain
\[
\mathbb{P}(M_t \leq y)=
\mathbb{P}(W_t\leq y, M_t \leq y) = 
 \Phi\left( \frac{y}{\sqrt{t}} \right) - \Phi\left( \frac{y-2y}{\sqrt{t}}\right) =2 \Phi\left( \frac{y}{\sqrt{t}} \right)  -1.\]

Due to the continuity of $\Phi$, we obtain
\[F_{WM}(x,y) =\Phi(x) - \Phi(x-2y).\]  
  
%\noindent
Case $0 \leq y < x$.

We apply the previous case and the inequality
\[ \forall t>0  \;\;\; W_t \leq M_t.\]

This leads to
\begin{eqnarray}
\mathbb{P}(W_t\leq x, M_t \leq y) &=& \mathbb{P}(W_t\leq y, M_t \leq y)\\
&=& \Phi\left( \frac{y}{\sqrt{t}} \right) - \Phi\left( \frac{y-2y}{\sqrt{t}}\right)= 2 \Phi\left( \frac{y}{\sqrt{t}} \right) -1. \nonumber
\end{eqnarray}
\end{proof}

\begin{Theorem}
The joint cumulative distribution function $F_{W,m,t}(x,z)$ of $(W_t,m_t)$, where $t>0$, is of the form
$$ F_{W,m,t}(x,z) =F_{W,m}\left( \frac{x}{\sqrt{t}},\frac{z}{\sqrt{t}}\right), \;\;\; t >0,$$
where
$$F_{W,m}(x,z)= \left\{
\begin{array}{ccc}
2 \Phi(z) -\Phi(2z-x) & \mbox{for } & z \leq min(0,x) ,\\
\Phi(x) & \mbox{for } &  z > \min(0,x).
\end{array}
\right. $$

Furthermore,
\[\mathbb{P}(m_t \leq z) = \min\left( 1, 2 \Phi\left(\frac{z}{\sqrt{t}} \right) \right).\]
\end{Theorem}

%\noindent
\begin{proof}
Case $x=z$.
We apply Formula (\ref{E:Wm}).
\begin{eqnarray}
\mathbb{P}(W_t\leq z, m_t \leq z) &=& 
\mathbb{P}(m_t \leq z) - \mathbb{P}(W_t \geq z, m_t \leq z) \nonumber \\
&=&  \mathbb{P}(m_t \leq z)-      \Phi\left(\frac{2z-z}{\sqrt{t}}\right). \nonumber
\end{eqnarray}

Since the event $\{W_t\leq z \}$ is contained in the event $\{m_t \leq z \}$,
\[ \mathbb{P}(m_t \leq z) = \Phi\left(\frac{z}{\sqrt{t}}\right) + \mathbb{P}(W_t\leq z, m_t \leq z)= 2 \Phi\left(\frac{z}{\sqrt{t}}\right).\]

Case $z \leq min(0,x)$.

We apply Formula (\ref{E:Wm}).
\begin{eqnarray}
\mathbb{P}(W_t\leq x, m_t \leq z) &=& 
\mathbb{P}(m_t \leq z) - \mathbb{P}(W_t \geq x, m_t \leq z) \\
&=&  2 \Phi\left(\frac{z}{\sqrt{t}}\right)-      \Phi\left(\frac{2z-x}{\sqrt{t}}\right). \nonumber
\end{eqnarray}

%\noindent
Case $z > x$.

The event $\{W_t\leq x \}$ is contained in the event $\{m_t \leq z \}$.
Indeed, for every $\omega$ from the first set,
\[ m_t(\omega) \leq W_t(\omega) \leq x \leq z.\]

Thus,
\[\mathbb{P}(W_t\leq x, m_t \leq z)=\mathbb{P}(W_t\leq x)= \Phi\left( \frac{x}{\sqrt{t}} \right).\]

%\noindent
Case $z > 0$.

Since $m_t \leq 0$, 
\[ \{\omega : m_t(\omega) \leq z\} = \Omega .\] 

Thus,
\[\mathbb{P}(W_t\leq x, m_t \leq z)=\mathbb{P}(W_t\leq x)= \Phi\left( \frac{x}{\sqrt{t}} \right).\]
\end{proof}

Note that, since the distributions of running maxima and minima are closely related with stopping times
\[ \mathbb{P} (M_t \leq y ) = 1 - \mathbb{P}(T_y \leq t),\]
\[ \mathbb{P} (m_t \leq z ) = \mathbb{P}(T_z \leq t),\]
they might be drawn from the formulas describing exit times (from the half-line) provided, for example, in~\cite{Schuss_2010} (Theorem 4.4.5).

\subsection{The Two-Barrier Stopping Time}\label{S:2barrier}

In this section, we apply the reflection principle for the two-barrier stopping time
\[\tau = \min (T_y,T_z), \;\; z<0<y,\]
and for the set $A$, which are equal, respectively, to 
\[ \{\omega: W_\tau(\omega) =y, \tau(\omega) \leq t \} \;\; \mbox{or} \;\; \{\omega: W_\tau(\omega) =z, \tau(\omega) \leq t \}.\]

Note that the range of $W_\tau$ consists only of two points $y$ and $z$.

\begin{Lemma}\label{L:Wleqx}
{For} 
 any $x \leq y$
\[ \mathbb{P} \left( W_t \leq x, \;\; W_\tau=y,\;\; \tau \leq t \right) =
\sum_{k=1}^\infty \left( \Phi\left(\frac{x + 2(k-1)z - 2ky}{\sqrt{t}} \right) - \Phi\left(\frac{x + 2kz - 2ky}{\sqrt{t}} \right) \right).\]
\end{Lemma}

\begin{proof}
Applying twice the reflection principle, we obtain %\vspace{-9pt}
%\begin{adjustwidth}{-\extralength}{0cm}
%\centering %% If there is a figure in wide page, please release command \centering
\begin{eqnarray}
\mathbb{P} \left( W_t \leq x, \;\; W_\tau=y,\;\; \tau \leq t \right) &=&
\mathbb{P} \left( 2 W_\tau -W_t \leq x, \;\; W_\tau=y,\;\; \tau \leq t \right)\\
&=& \mathbb{P} \left( 2 y -W_t \leq x, \;\; W_\tau=y,\;\; \tau \leq t \right) \nonumber \\
&=& \mathbb{P} \left( W_t  \geq 2y- x, \;\; W_\tau=y,\;\; \tau \leq t \right) \nonumber\\
&=& \mathbb{P} \left( W_t  \geq 2y- x, \;\; \tau \leq t \right)
- \mathbb{P} \left( W_t  \geq 2y- x, \;\; W_\tau=z,\;\; \tau \leq t \right) \nonumber\\
&=& \mathbb{P} \left( W_t  \geq 2y- x, \;\; \tau \leq t \right)
- \mathbb{P} \left( 2W_\tau - W_t  \geq 2y- x, \;\; W_\tau=z,\;\; \tau \leq t \right) \nonumber\\
&=& \mathbb{P} \left( W_t  \geq 2y- x, \;\; \tau \leq t \right)
- \mathbb{P} \left( W_t  \leq 2z -2y+ x, \;\; W_\tau=z,\;\; \tau \leq t \right) \nonumber\\
&=& \mathbb{P} \left( W_t  \geq 2y- x, \;\; \tau \leq t \right)
- \mathbb{P} \left( W_t  \leq 2z -2y+ x, \;\; \tau \leq t \right) \nonumber \\
&&+ \, \mathbb{P} \left( W_t  \leq 2z -2y+ x, \;\; W_\tau=y,\;\; \tau \leq t \right) .\nonumber
\end{eqnarray}
%%\end{adjustwidth}

Since $x \leq y$, the event $ W_t  \geq 2y- x $ is contained in the event $\tau \leq t$, and
\[ W_t(\omega) \geq 2y-x \geq y \Longrightarrow M_t(\omega) \geq y \Longrightarrow T_y(\omega) \leq t \Longrightarrow \tau(\omega) =\min(T_y(\omega),T_z(\omega))\leq t.\]

Therefore,
\[\mathbb{P} \left( W_t  \geq 2y- x, \;\; \tau \leq t \right)=\mathbb{P} \left( W_t  \geq 2y- x \right) = \Phi_t(x-2y),\]
where $\Phi_t$ denotes the distribution function of an $N(0,t)$ random variable:
\[\Phi_t(x)=\Phi\left( \frac{x}{\sqrt{t}}\right).\]

Similarly, since $z<0<y$ and $x \leq y$, the event $ W_t  \leq 2z -2y+ x $ is contained in the event $\tau \leq t$, so
%\begin{adjustwidth}{-\extralength}{0cm}
%\centering %% If there is a figure in wide page, please release command \centering
\[ W_t(\omega) \leq 2z -2y+ x \leq z  \Longrightarrow m_t(\omega) \leq z \Longrightarrow T_z(\omega) \leq t \Longrightarrow \tau(\omega) =\min(T_y(\omega),T_z(\omega))\leq t.\]
%%\end{adjustwidth}

Therefore,
\[\mathbb{P} \left( W_t  \leq 2z -2y+ x, \;\; \tau \leq t \right)=\mathbb{P} \left( W_t  \leq 2z -2y+ x \right) = \Phi_t(x+ 2z-2y).\]

In such a way, we obtain the recurrence for $x\leq y$:
%\begin{adjustwidth}{-\extralength}{0cm}
%\centering %% If there is a figure in wide page, please release command \centering
\begin{equation}
\mathbb{P} \left( W_t \leq x, \; W_\tau=y,\; \tau \leq t \right) =
\mathbb{P} \left( W_t  \leq 2z -2y+ x, \;W_\tau=y,\; \tau \leq t \right)
+ \Phi_t(x-2y) - \Phi_t(x+ 2z-2y).
\end{equation}
%%\end{adjustwidth}

Since $2z-2y <0$, we may repeat it. After $n$ repetitions, we obtain
%\begin{adjustwidth}{-\extralength}{0cm}
%\centering %% If there is a figure in wide page, please release command \centering
\begin{eqnarray}
\mathbb{P} \left( W_t \leq x, \; W_\tau=y,\; \tau \leq t \right) &=&
\mathbb{P} \left( W_t  \leq n(2z -2y)+ x, \;W_\tau=y,\; \tau \leq t \right)\\
&&+\; \sum_{k=1}^n \left( \Phi_t(x+ 2(k-1)z -2ky) - \Phi_t(x+ 2kz-2ky) \right). \nonumber
\end{eqnarray}
%%\end{adjustwidth}

Since 
\[\mathbb{P} \left( W_t  \leq n(2z -2y)+ x, \;W_\tau=y,\; \tau \leq t \right) \leq \mathbb{P} \left( W_t  \leq n(2z -2y)+ x \right) \stackrel{n\rightarrow \infty}{\longrightarrow } 0,\]
and passing to the limit $n\rightarrow \infty$, we obtain
\begin{equation}
\mathbb{P} \left( W_t \leq x, \; W_\tau=y,\; \tau \leq t \right) =
\sum_{k=1}^\infty \left( \Phi_t(x+ 2(k-1)z -2ky) - \Phi_t(x+ 2kz-2ky) \right).
\end{equation}
\end{proof}

\begin{Lemma}\label{L:Wgeqx}
For any $x \geq z$
\[ \mathbb{P} \left( W_t \geq x, \;\; W_\tau=z,\;\; \tau \leq t \right) =
%\sum_{k=1}^\infty \left( \Phi\left(\frac{-x + 2kz - 2(k-1)y}{\sqrt{t}} \right) - \Phi\left(\frac{-x + 2kz - 2ky}{\sqrt{t}} \right) \right).\]
\sum_{k=1}^\infty \left( \Phi\left(\frac{x - 2kz + 2ky}{\sqrt{t}} \right)  - \Phi\left(\frac{x - 2kz + 2(k-1)y}{\sqrt{t}} \right) \right).\]
\end{Lemma}

\begin{proof}
We base our proof on the symmetry of the Wiener process $W_t$. 
We put $\tilde{W}_t=-W_t$. Then, we have, for $z<0<y$,
\begin{eqnarray*}
\tilde{M}_t &=& \sup_{0 \le s \le t}\tilde{W}_s = - \inf_{0 \le s \le t}W_s = -m_t,\\
\tilde{m}_t &=& \inf_{0 \le s \le t}\tilde{W}_s=- \sup_{0 \le s \le t}W_s =-M_t,\\
\tilde{T}_{-z}&=&\inf \{t> 0: \tilde{W}_t = -z\}=T_z,\\
\tilde{T}_{-y}&=&\inf \{t> 0: \tilde{W}_t = -y\}=T_y,\\
\tilde{\tau}&=& \min(\tilde{T}_{-z},\tilde{T}_{-y})=\min(T_z,T_y)= \tau .
\end{eqnarray*}

Therefore,
\[ \mathbb{P} \left( W_t \geq x, \;\; W_\tau=z,\;\; \tau \leq t \right) =\mathbb{P} \left( \tilde{W}_t \leq -x, \;\; \tilde{W}_{\tilde{\tau}}=-z,\;\; \tilde{\tau} \leq t \right) .\]

Since $\tilde{W}_t$ is as well a Wiener process, we apply Lemma \ref{L:Wleqx} and obtain 
%\begin{adjustwidth}{-\extralength}{0cm}
%\centering %% If there is a figure in wide page, please release command \centering
\begin{eqnarray*}
 \mathbb{P} \left( \tilde{W}_t \leq -x, \;\; \tilde{W}_{\tilde{\tau}}=-z,\;\; \tilde{\tau} \leq t \right) &=&
\sum_{k=1}^\infty \left( \Phi\left(\frac{-x + 2kz - 2(k-1)y}{\sqrt{t}} \right) - \Phi\left(\frac{-x + 2kz - 2ky}{\sqrt{t}} \right) \right)  \\
&=&\sum_{k=1}^\infty \left( \Phi\left(\frac{x - 2kz + 2ky}{\sqrt{t}} \right)  - \Phi\left(\frac{x - 2kz + 2(k-1)y}{\sqrt{t}} \right) \right). \;\;\;\;
\end{eqnarray*}
%%\end{adjustwidth}
\end{proof}

Together, the above lemmas  give us

\begin{Corollary}\label{C:W_tau}
{For} 
 $z<0<y$ and $z \leq x \leq y$, %\vspace{-12pt}
%\begin{adjustwidth}{-\extralength}{0cm}
%\centering %% If there is a figure in wide page, please release command \centering
\begin{eqnarray}
\mathbb{P} \left( W_t \leq x, \; \tau \leq t \right) &=&
\mathbb{P} \left( W_\tau =z, \; \tau \leq t \right)\\
&&+ \sum_{k=-\infty}^\infty \left( \Phi\left(\frac{x + 2(k-1)z - 2ky}{\sqrt{t}} \right) - \Phi\left(\frac{x + 2kz - 2ky}{\sqrt{t}} \right) \right) %\nonumber \\
+ \Phi(x) \nonumber
\end{eqnarray}
%%\end{adjustwidth}
%-\sum_{k=1}^\infty \left( \Phi\left(\frac{-x + 2kz - 2(k-1)y}{\sqrt{t}} \right) - \Phi\left(\frac{-x + 2kz - 2ky}{\sqrt{t}} \right) \right) , \nonumber\\
where 
%\begin{adjustwidth}{-\extralength}{0cm}
%\centering %% If there is a figure in wide page, please release command \centering
\begin{eqnarray}
\mathbb{P} \left( W_\tau =z, \; \tau \leq t \right) &=&
2\sum_{k=1}^\infty \left( \Phi\left(\frac{(2k-1)z - 2(k-1)y}{\sqrt{t}} \right) - \Phi\left(\frac{(2k-1)z - 2ky}{\sqrt{t}} \right) \right) .
\end{eqnarray}
%%\end{adjustwidth}

Furthermore, %\vspace{-12pt}
%\begin{adjustwidth}{-\extralength}{0cm}
%\centering %% If there is a figure in wide page, please release command \centering
\begin{eqnarray}
\mathbb{P} \left( W_\tau =y, \; \tau \leq t \right) &=&
2\sum_{k=1}^\infty \left( \Phi\left(\frac{2(k-1)z - (2k-1)y}{\sqrt{t}} \right) - \Phi\left(\frac{2kz - (2k-1)y}{\sqrt{t}} \right) \right)
\end{eqnarray}
%%\end{adjustwidth}
and %\vspace{-12pt}
%\begin{adjustwidth}{-\extralength}{0cm}
%\centering %% If there is a figure in wide page, please release command \centering
\begin{eqnarray}
\mathbb{P} \left( \tau \leq t \right) &=&
2\sum_{k=1}^\infty \left( \Phi\left(\frac{(2k-1)z - 2(k-1)y}{\sqrt{t}} \right) - \Phi\left(\frac{(2k-1)z - 2ky}{\sqrt{t}} \right) \right) \label{E:tau}\\
&&+2\sum_{k=1}^\infty \left( \Phi\left(\frac{2(k-1)z - (2k-1)y}{\sqrt{t}} \right) - \Phi\left(\frac{2kz - (2k-1)y}{\sqrt{t}} \right) \right). \nonumber
\end{eqnarray}
%%\end{adjustwidth}

\end{Corollary}

\begin{proof}
We observe that
\begin{eqnarray*}
\mathbb{P} \left( W_t \leq x, \; \tau \leq t \right) &=&
\mathbb{P} \left( W_t \leq x, \; W_\tau=y, \; \tau \leq t \right)  + \mathbb{P} \left( W_t \leq x, \; W_\tau=z, \; \tau \leq t \right) \\
&=&
\mathbb{P} \left( W_t \leq x, \; W_\tau=y, \; \tau \leq t \right)\\
&& +
\mathbb{P} \left( W_\tau =z, \; \tau \leq t \right) - \mathbb{P} \left( W_t \geq x, \; W_\tau=z, \; \tau \leq t \right).
\end{eqnarray*}

Thus, applying Lemmas \ref{L:Wleqx} and \ref{L:Wgeqx}, we obtain
\begin{eqnarray}
\mathbb{P} \left( W_t \leq x, \; \tau \leq t \right) &=&
\mathbb{P} \left( W_\tau =z, \; \tau \leq t \right)\\
&&+ \sum_{k=1}^\infty \left( \Phi_t(x + 2(k-1)z - 2ky ) - \Phi_t(x + 2kz - 2ky ) \right)  \nonumber \\
&&+
\sum_{k=1}^\infty \left( \Phi_t(x - 2kz + 2(k-1)y) - \Phi_t(x - 2kz + 2ky) \right) . \nonumber \\
&=& 
\mathbb{P} \left( W_\tau =z, \; \tau \leq t \right)  + \Phi(x) \nonumber\\
&&+ \sum_{k=-\infty}^\infty \left( \Phi_t(x + 2(k-1)z - 2ky ) - \Phi_t(x + 2kz - 2ky ) \right)  \nonumber
\end{eqnarray}

Since the event $W_t \leq z$ implies the event $\tau \leq t$, we have
\[ \mathbb{P} \left( W_t \leq z, \; \tau \leq t \right)= \mathbb{P} \left( W_t \leq z \right) = \Phi_t(z).\]

Therefore, substituting $x=z$ in the above formula, we obtain
\begin{eqnarray*}
\mathbb{P} \left( W_\tau =z, \; \tau \leq t \right) &=&
\mathbb{P} \left( W_t \leq z, \; \tau \leq t \right) -\Phi(z)\\
&&- \sum_{k=-\infty}^\infty \left( \Phi_t(z + 2(k-1)z - 2ky ) - \Phi_t(z + 2kz - 2ky ) \right)  \\
%&& + \sum_{k=1}^\infty \left( \Phi_t(-z + 2kz - 2(k-1)y) - \Phi_t(-z + 2kz - 2ky) \right)\\
&=& \Phi_t(z) - 2 \sum_{k=1}^\infty \Phi_t( (2k-1)z - 2ky )\\
&&+\sum_{k=2}^\infty \Phi_t( (2k-1)z - 2(k-1)y )
+\sum_{k=1}^\infty \Phi_t( (2k-1)z - 2(k-1)y )\\
%&=& \Phi_t(z) - 2 \sum_{k=1}^\infty \Phi_t( (2k+1)z - 2ky ) + 2 \sum_{k=1}^\infty \Phi_t( (2k-1)z - 2ky ) - \Phi_t(z)\\
&=& 2 \sum_{k=1}^\infty \left( \Phi_t( (2k-1)z - 2(k-1)y )  - \Phi_t( (2k-1)z - 2ky ) \right).
\end{eqnarray*}

This concludes the proof of the first two formulas of the corollary.
The third one follows from the symmetry of the Wiener process. We apply the notation from the proof of Lemma \ref{L:Wgeqx}.

Since $\tilde{W}_t= -W_t$ is a Wiener process,%\vspace{-12pt}
%\begin{adjustwidth}{-\extralength}{0cm}
%\centering %% If there is a figure in wide page, please release command \centering
\begin{eqnarray*}
\mathbb{P} \left( W_\tau =y, \; \tau \leq t \right) &=&
\mathbb{P} \left( \tilde{W}_{\tilde{\tau}} =-y, \; \tilde{\tau} \leq t \right)\\
&=& 2 \sum_{k=1}^\infty \left( \Phi_t( (2k-1)(-y) - 2(k-1)(-z) )  - \Phi_t( (2k-1)(-y) - 2k(-z) ) \right)\\
&=& 2 \sum_{k=1}^\infty \left( \Phi_t(  2(k-1)z  - (2k-1)y)  - \Phi_t( 2kz - (2k-1)y  ) \right).
\end{eqnarray*}
%%\end{adjustwidth}

Since
\[ \mathbb{P} \left(\tau \leq t \right) =\mathbb{P} \left( W_\tau =y, \; \tau \leq t \right) + \mathbb{P} \left( W_\tau =z, \; \tau \leq t \right), \]
the last equality of the corollary is a consequence of the two previous ones.
\end{proof}

Note that, since the Formula (\ref{E:tau}) is describing an exit time from a segment, it might be obtained
as a solution of a partial differential equation (see~\cite{Schuss_2010} Theorem 4.4.5).

\subsection{Joint Cumulative Distribution}\label{S:Proof_1}

\begin{Theorem}\label{T:F}
The joint cumulative distribution function $F_t(x,y,z)$ of $(W_t, M_t,m_t)$, where $t>0$, is of the~form
$$ F_t(x,y,z) =F\left( \frac{x}{\sqrt{t}},\frac{y}{\sqrt{t}}, \frac{z}{\sqrt{t}} \right), \;\;\; t >0,$$
where
%\begin{adjustwidth}{-\extralength}{0cm}
%\centering %% If there is a figure in wide page, please release command \centering
\begin{equation}\label{E:F}
F(x,y,z)= \left\{
\begin{array}{cl}
0 & \mbox{for }   y \leq 0,\\
\Phi(x)-\Phi(x-2y) & \mbox{for }   z \geq 0, \;  0 < y,  \;  x  < y,\\
 2\Phi(y)-1& \mbox{for }   z \geq 0, \;  0 < y \leq x, \\
 \Phi(x)-\Phi(x-2y) & \mbox{for }   x <z <0 <y,  \\
 2\Psi(z,2y,2(y-z)) - \Psi(x+2z-2y,2y,2(y-z)) &\\
-\Psi(-x+2z,2y,2(y-z)) & \mbox{for }   z <0 <y, \; z \leq x \leq y,\\
2\Psi(z,2y,2(y-z)) - 2\Psi(2z-y,2y,2(y-z))  &\mbox{for }   z <0 <y < x ,
\end{array}
\right. \end{equation}
%%\end{adjustwidth}
where
\begin{equation}\label{E:Psi1}
\Psi(q,r,s) = \sum_{k=0}^\infty (\Phi(q-ks) - \Phi(q-r-ks)) 
\end{equation}
and $\Phi$ denotes the distribution function of the standard normal probability law N(0, 1). 
\end{Theorem}

Note that the distribution of $(W_t,M_t,m_t)$ is absolutely continuous with respect to the Lebesgue measure (Section \ref{S:Properties_1}) and is supported on the set
\[ \{(x,y,z) \in \mathbb{R}^3:\; z\leq 0 \leq y, \; z \leq x \leq y\}.\]

Furthermore, 
the marginal distributions of $F$ are given by the following formulas:
%\begin{adjustwidth}{-\extralength}{0cm}
%\centering %% If there is a figure in wide page, please release command \centering
\begin{eqnarray}
F(x,+\infty,+\infty) &=& \mathbb{P}(W_1\leq x) = \Phi(x),\\
F(+\infty,y,+\infty) &=& \mathbb{P}(M_1 \leq y)= \left\{\begin{array}{cl}
0 & \mbox{if } \; y \leq 0,\\
2\Phi(y)-1 & \mbox{if } \; y > 0,
\end{array} \right.\\
F(+\infty,+\infty,z) &=& \mathbb{P}(m_1 \leq z)= \left\{\begin{array}{cl}
2 \Phi(z) & \mbox{if } \; z< 0,\\
1 & \mbox{if } \; z \geq  0.
\end{array} \right.\\
F(x,y,+\infty) &=& \mathbb{P}(W_1\leq x, M_1 \leq y)= \left\{\begin{array}{cl}
0 & \mbox{if } \; y \leq 0,\\
2\Phi(y)-1 & \mbox{if } \; 0<y  \leq x,\\
\Phi(x) - \Phi(x-2y) & \mbox{if } \; \max(0,x) <y,
\end{array} \right.\\
F(x,+\infty, z) &=& \mathbb{P}(W_1\leq x, m_1 \leq z)= \left\{\begin{array}{cl}
2\Phi(z) - \Phi(2z-x) & \mbox{if } \; z \leq \min(0,x),\\
\Phi(x) & \mbox{if } \; x<z  \leq 0  \vee 0<z,
\end{array} \right.\\
F(+\infty, y, z) &=& \mathbb{P}(M_1\leq y, m_1 \leq z)= \left\{\begin{array}{cl}
0 & \mbox{if } \;  y \leq 0,\\
2\Psi(z,2y,2(y-z)) &\\
- 2\Psi(2z-y,2y,2(y-z))& \mbox{if } \; z  \leq 0  < y,\\
2 \Phi(y)-1 & \mbox{if } \; 0 <z, 0<y.
\end{array} \right.
\end{eqnarray}
%%\end{adjustwidth}

We prove Theorem \ref{T:F} case by case.

%\noindent
\begin{proof}
{Case} 
 $y \leq 0$.

Since $W_t$ is almost surely continuous, and $W_0=0$, $M_t$ is almost surely positive:
\[F_t(x,y,z) =\mathbb{P}( W_t \leq x, \;  M_t \leq y, \; m_t \leq z) \leq  \mathbb{P}(M_t\leq y) =0.\]

%\noindent
Case $z \geq 0$ and $y>0$.

Since $m_t$ is almost surely negative,
\[F_t(x,y,z) =\mathbb{P}( W_t \leq x, \; M_t \leq y, \; m_t \leq z) =\mathbb{P}( W_t \leq x, \; M_t \leq y).\]

Following Theorem \ref{T:F_WM}, we obtain
\[F_t(x,y,z) =
\left\{
\begin{array}{ccc}
2 \Phi_t(y) -1 & \mbox{for } &  y < x,\\
\Phi_t(x) - \Phi_t(x-2y) & \mbox{for } &  x \leq y .
\end{array}
\right. \]

%\noindent
Case $z < 0 < y$ and  $x \leq z $.

Since the event $W_t \leq x \leq z$ implies the event $m_t \leq z$,
 \[F_t(x,y,z) =\mathbb{P}( W_t \leq x, \; M_t \leq y, \; m_t \leq z) =\mathbb{P}( W_t \leq x, \; M_t \leq y).\]
 
 As above, following Theorem \ref{T:F_WM}, we obtain
\[F_t(x,y,z) = \Phi_t(x) - \Phi_t(x-2y) .\]

%\noindent
Case $z < 0 < y$ and  $z<x \leq y$.

We observe that
\begin{eqnarray*}
F_t(x,y,z) &=&\mathbb{P}( W_t \leq x, \; M_t \leq y, \; m_t \leq z) \\
&=& \mathbb{P}( W_t \leq x, \;( M_t \geq y \vee m_t \leq z))\\
&& - \mathbb{P}( W_t \leq x, \; M_t \geq y)\\
&=& \mathbb{P}( W_t \leq x, \; \tau \leq t)\\
&& - \mathbb{P}( W_t \leq x)  +\mathbb{P}( W_t \leq x, \; M_t \leq y) 
\end{eqnarray*}
where $\tau$ is the two-barrier stopping time discussed in Section \ref{S:2barrier}.
Following Corollary \ref{C:W_tau} and Theorem \ref{T:F_WM}, we obtain
\begin{eqnarray*}
F_t(x,y,z) &=&
2 \sum_{k=1}^\infty \left( \Phi_t( (2k-1)z - 2ky )  - \Phi_t( (2k+1)z - 2ky ) \right)\\
&&+ \sum_{k=1}^\infty \left( \Phi_t(x + 2(k-1)z - 2ky ) - \Phi_t(x + 2kz - 2ky ) \right)  \\
&&- \sum_{k=1}^\infty \left( \Phi_t(-x + 2kz - 2(k-1)y) - \Phi_t(-x + 2kz - 2ky) \right) \\
&& -\Phi_t(x) + \Phi_t(x) - \Phi_t(x-2y)\\
&=& 2 \sum_{k=1}^\infty \left( \Phi_t( (2k-1)z - 2ky )  - \Phi_t( (2k+1)z - 2ky ) \right)\\
&&+ \sum_{k=1}^\infty \left( \Phi_t(x + 2kz - 2(k+1)y )  - \Phi_t(x + 2kz - 2ky ) \right)  \\
&&- \sum_{k=1}^\infty \left( \Phi_t(-x + 2kz - 2(k-1)y) - \Phi_t(-x + 2kz - 2ky) \right) .
\end{eqnarray*}
 
% \noindent
Case $z < 0 < y$ and  $y \leq x$.

Since $W_t \leq M_t$,
\begin{eqnarray*}
F_t(x,y,z) &=&\mathbb{P}( W_t \leq x, \; M_t \leq y, \; m_t \leq z)=\mathbb{P}( M_t \leq y, \; m_t \leq z)\\
&=& \mathbb{P}( M_t \geq y \vee  m_t \leq z) + \mathbb{P}(M_t \leq y) -1\\
&=& \mathbb{P}( \tau \leq t) + \mathbb{P}(M_t \leq y) -1\\
&=& 2 \sum_{k=1}^\infty \left( \Phi_t( (2k-1)z - 2(k-1)y )  - \Phi_t( (2k-1)z - 2ky ) \right)\\
&& + 2 \sum_{k=1}^\infty \left( \Phi_t(  2(k-1)z  - (2k-1)y)  - \Phi_t( 2kz - (2k-1)y  ) \right)\\
&& - 2 \Phi(-y) \\
&=& 2 \sum_{k=1}^\infty \left( \Phi_t( (2k-1)z - 2(k-1)y )  - \Phi_t( (2k-1)z - 2ky ) \right)\\
&& - 2 \sum_{k=1}^\infty \left( \Phi_t(  2kz  - (2k-1)y)  - \Phi_t( 2kz - (2k+1)y  ) \right)
\end{eqnarray*}
\end{proof}

{
As a corollary of Theorem \ref{T:F}, we prove the 80-year-old characterization of the Wiener process called the ``L{\'e}vy Triple Law'' (\cite{Schilling_2014} Theorem 6.18), which deals with the trajectories lying between the upper barrier $u$ and lower barrier $l$.

%\pagebreak
\begin{Corollary}\label{C:LTL} 
For $l<0<u$, $l \leq x \leq u$, and $t>0$,
\begin{equation}\label{E:LTL}
 \frac{\partial}{\partial x} \mathbb{P}(W_t \leq x , M_t \leq u , m_t \geq l) = \sum_{k=-\infty}^{+\infty} \left( \varphi_t(x +2k(u-l)) - \varphi_t(x-2l -2k(u-l) \right) ,
 \end{equation}
where $\varphi_t$ denotes the density of $N(0,t)$.
\end{Corollary}

%\noindent
\begin{proof}
We observe that
\begin{eqnarray}
&& \mathbb{P}(W_t \leq x , M_t \leq u , m_t \geq l)= F_t(x,u,0)-F_t(x,u,l)\\
 && \;\;\;\;  = \Phi_t(x)-\Phi_t(x-2u) \nonumber \\
&&\;\;\;\;\;\; - \sum_{k=0}^{+\infty} \left( \Phi_t(l -2k(u-l)) - \Phi_t(l - 2u -2k(u-l)) \right) \nonumber \\
&&\;\;\;\;\;\; + \sum_{k=0}^{+\infty} \left( \Phi_t(x+2l-2u -2k(u-l)) - \Phi_t(x+2l-2u - 2u -2k(u-l)) \right) \nonumber \\
&&\;\;\;\;\;\; + \sum_{k=0}^{+\infty} \left( \Phi_t(-x+2l -2k(u-l)) - \Phi_t(-x+2l - 2u -2k(u-l)) \right). \nonumber 
\end{eqnarray}

Hence,
\begin{eqnarray}
 && \frac{\partial}{\partial x} \mathbb{P}(W_t \leq x , M_t \leq u , m_t \geq l)\\
  && \;\;\;\;  = \varphi_t(x)-\varphi_t(x-2u) \nonumber \\
&&\;\;\;\;\;\; + \sum_{k=0}^{+\infty} \left( \varphi_t(x+2l-2u -2k(u-l)) - \varphi_t(x+2l-2u - 2u -2k(u-l)) \right) \nonumber \\
&&\;\;\;\;\;\; - \sum_{k=0}^{+\infty} \left( \varphi_t(-x+2l -2k(u-l)) - \varphi_t(-x+2l - 2u -2k(u-l)) \right). \nonumber \\
&& \;\;\;\;  = \varphi_t(x)-\varphi_t(x-2u) \nonumber \\
&&\;\;\;\;\;\; + \sum_{k=0}^{+\infty} \left( \varphi_t(x -2(k+1)(u-l)) - \varphi_t(x - 2u -2(k+1)(u-l)) \right) \nonumber \\
&&\;\;\;\;\;\; - \sum_{k=0}^{+\infty} \left( \varphi_t(x-2u +2(k+1)(u-l)) - \varphi_t(x + 2(k+1) (u-l)) \right). \nonumber \\
&&\;\;\;\;\;\; = \sum_{k=-\infty}^{+\infty}  \varphi_t(x -2k(u-l)) - \sum_{k=-\infty}^{+\infty}\varphi_t(x - 2u -2k(u-l))  \nonumber \\
&&\;\;\;\;\;\; = \sum_{k=-\infty}^{+\infty}  \varphi_t(x -2k(u-l)) - \sum_{k=-\infty}^{+\infty}\varphi_t(x - 2l -2(k+1)(u-l)) . \nonumber 
\end{eqnarray}

{After shifting by one the parameter $k$ in the second sum, we obtain the required \mbox{Formula~(\ref{E:LTL})}.}

\end{proof}
}

%\pagebreak

\subsection{Properties of the Cumulative Distribution Function $F$ and the Marginal Distributions}\label{S:Properties_1}

We start with the absolute continuity of the marginal distributions of random pairs $(W_1,M_1)$ and $(W_1,m_1)$:
\[F_{W,M}(x,y)=F(x,y,\infty) \;\; \mbox{and } \;\; F_{W,m}(x,z)=F(x,\infty, z).\]

\begin{Proposition}
The bivariate cumulative distributions $F_{W,M}$ and  $F_{W,m}$ are absolutely continuous with respect to the Lebesgue measure.
\end{Proposition}

\begin{proof}
We start with the function $F_{W,m}(x,z)$.
We show that a function
\[ g(\xi,\eta) = -2 (2\eta- \xi) \varphi(2\eta - \xi) \jedynka_{\eta \leq 0} \jedynka_{\eta \leq \xi}\]
is a corresponding density. Let $I(x,z)$ denote an integral of $g$ over the set $(-\infty,x]\times (-\infty,z]$. 
\begin{eqnarray}
I(x,y) &=& \int_{-\infty}^z \int_{-\infty}^x g(\xi,\eta) d\eta d \xi\\
&=& \int_{-\infty}^{\min(z,x,0)} \int_{\eta}^x (-2) (2\eta- \xi) \varphi(2\eta - \xi) d\eta d \xi \nonumber \\
&=& \int_{-\infty}^{\min(z,x,0)} 2 \bigg( \varphi(\eta) - \varphi(2\eta - x) \bigg) d\xi \nonumber \\
&=& 2 \Phi(\min(z,x,0)) - \Phi(2 \min(z,x,0)-x) \nonumber \\
&=& \left\{\begin{array}{cl}
2 \Phi(z) - \Phi(2z-x) & \mbox{when } \; \min(z,x,0)=z,\\
2 \Phi(x) - \Phi(2x-x) = \Phi(x)& \mbox{when } \; \min(z,x,0)=x,\\
1-\Phi(0-x)=\Phi(x) & \mbox{when } \; \min(z,x,0)=0.
\end{array}  \right. \nonumber
\end{eqnarray}

Since all points $(x,z)$, $\, I(x,z)$, and $F_{W,m}(x,z)$ are equal each other, we conclude that $F_{W,m}$ is absolutely continuous.

Next, since $\tilde{W}_t= -W_t$ is a Wiener process, the pairs $(W_1,M_t)$ and $(-W_1,-m_1)$ have the same distribution
\[ (W_1,M_t) \stackrel{d}{=} (-W_1,-m_1).\]

Therefore, $F_{W,M}$ is also absolutely continuous.
\end{proof}

Now, we are in position to show the absolute continuity of the cumulative distribution function $F$ of the random triple $(W_1,M_1,m_1)$.

\begin{Proposition}
The function $F$ is absolutely  continuous with respect to the Lebesgue measure. Furthermore,
\[\mathbb{P}\bigg( (W_1,M_1,m_1) \in \{(x,y,z) \in \mathbb{R}^3 : z<0 < y, z<x<y\} \bigg)=1 .\]
\end{Proposition}

\begin{proof}
Since {the definition of running maxima and minima implies that}
\[ m_1 \leq W_1 \leq M_1, \; \mbox{and } \; m_1 \leq 0 \leq M_1,\]
we have
\[\mathbb{P}\bigg( (W_1,M_1,m_1) \in \{(x,y,z) \in \mathbb{R}^3 : z\leq 0 \leq y, z \leq x \leq y\}\bigg) =1 .\]

We show that the probability that the random triple $ (W_1,M_1,m_1) $ belongs to the boundary and is equal to 0.

Since the distribution functions of $M_1$ and $m_1$ are continuous, we have
\[ \mathbb{P}(M_1=0) =0, \;\; \mathbb{P}(m_1=0) =0.\]

Since the distribution functions of the random pairs $(W_1,M_1)$ and $(W_1,m_1)$ are absolutely continuous, 
\[ \mathbb{P}(M_1=W_1) =0 = \mathbb{P}(m_1=W_1) .\]

To conclude the proof, it is enough to observe that the function $F$  is $C^\infty$ on the set
\[ \{(x,y,z) \in \mathbb{R}^3 : z<0 < y, z<x<y\}. \]

Hence, it is absolutely continuous with respect to the Lebesgue measure.
\end{proof}

%\pagebreak

\subsection{Proof of Theorem \ref{T:C}}\label{S:Proof_2}
We recall that a \emph{{copula}%%%: Please confirm if the italics is unnecessary and can be removed. please check whole text
} is the restriction to the unit $n$-cube $[0,1]^n$ of a distribution function  whose univariate
margins are uniformly distributed on $[0,1]$. Specifically, a function $C:[0,1]^n\rightarrow [0,1]$ is a copula if it 
has the following properties.
For every ${\it \bf u}=(u_1, \dots , u_n)$ and ${\it \bf v}=(v_1, \dots , v_n)$ such that
$0 \leq u_i \leq v_i \leq 1 $ for $i=1, \dots , n$, we have the following:
\begin{itemize}
\item[(C1) ] $(\exists i \;\; u_i=0 ) \Rightarrow C({\it \bf u})=0$; %i.e$.$ $C$ is \emph{grounded},
\item[(C2) ] $\forall j \in \{1, \dots , n\} \;\;\;( \forall i \neq j \;\; u_i=1 ) \Rightarrow C({\it \bf u})=u_j$;
\item[(C3) ] $C$ is \emph{$n$-nondecreasing}, that is, the $C$-volume $V_C ({\it \bf u},{\it \bf v})$ of any $n$-rectangle
with lower vertex ${\it \bf u}$ and upper vertex ${\it \bf v}$ is non-negative.
\end{itemize}

We recall that the $C$-volume %Please check meaning retained. 
is a signed sum of the values of $C$ at the vertices of the $n$-rectangle,%\vspace{-6pt}
\[V_C ({\it \bf u},{\it \bf v})= 
\left.C({\it \bf w})\right\vert_{w_1=u_1}^{v_1} \left. \dots \right\vert_{w_n=u_n}^{v_n}=
\sum_{j_1=1}^{2} \dots \sum_{j_n=1}^{2} (-1)^{j_1 + \dots + j_n} C(w_{1,j_1}, \dots , w_{n,j_n}),\]
where 
$w_{i,1}=u_i$, and $w_{i,2}=v_i$.

According to the celebrated \emph{Sklar's theorem}, the joint distribution function $F$ of any $n$-tuple 
 ${\it \bf X}=(X_1, \dots , X_n)$ of  random variables defined on the
probability space $(\Omega,\mathcal{F},\textbf{P})$ can be written as a composition of a copula $C$ 
and the univariate marginals $F_i$, i.e., for all ${\it \bf x}=(x_1, \dots , x_n)\in \mathbb{R}^n$, 
\[F({\it \bf x})=C(F_1(x_1),\dots , F_n(x_n)).\]

Moreover, if the distribution functions $F_k$ are continuous, then the copula $C$ is uniquely determined and can be described in terms of the quantile functions, i.e., the generalized inverses of $F_k$. Since $F_k(F_k^\leftarrow (u_k)=u_k$,  we obtain 
\[C(u_1, \dots , u_n) = C(F_1(F_1^\leftarrow (u_1)),\dots , F_n(F_n^\leftarrow (u_n)) )= F(F_1^\leftarrow (u_1), \dots , F_n^\leftarrow (u_n).\]

In our case where $n=3$ and for $u,v,w \in (0,1)$, the quantile functions of $W_1, M_1$, and $m_1$ are given by %\vspace{-6pt}
\begin{eqnarray}
F_1^\leftarrow (u) &=& \Phi^{-1}(u),\\
F_2^\leftarrow (v)&=& \Phi^{-1}\left(\frac{1+v}{2}\right),\\
F_3^\leftarrow (w) &=& \Phi^{-1}\left( \frac{w}{2} \right).
\end{eqnarray}

To conclude the proof of Theorem \ref{T:C}, we substitute the above quantile functions into Formula (\ref{E:F}).

%\pagebreak

\subsection{Properties of the Copula $C_{W,M,m}$ and Its Marginal Copulas}\label{S:Properties}

The copulas mentioned in Section \ref{S:Main} regarding the involution of the copula $C_{W,M,m}$ and its marginal copulas follows from the fact that %Please check meaning retained.
{$-W_t$ is a Wiener process as well. Hence, the triple $(-W_1,-m_1,-M_1)$ has the same distribution as the triple $(W_1,M_1,m_1)$.}
We put
\[ U= \Phi(W_1), \; V= 2 \Phi(M_1)-1, \; W= 2 \Phi(m_1) .\]

Obviously, $(U,V,W)$ are the generators of the copula 
$C_{W,M,m}$
\[ C_{W,M,m}(u,v,w) = \mathbb{P}(U\leq u, V \leq v, W \leq w).\]

Since %\vspace{-6pt}
\begin{eqnarray*}
\Phi(-W_1) &=&  1 - \Phi(W_1) = 1-U,\\
2 \Phi(-m_1) -1 &=& 2( 1 - \Phi(m_1)) -1 = 1 - W,\\
2 \Phi(-M_1) &=& 2 (1- \Phi(M_1) )  =1 -V,
\end{eqnarray*}
we apply the inclusion/exclusion principle and obtain
%\begin{adjustwidth}{-\extralength}{0cm}
%\centering %% If there is a figure in wide page, please release command \centering
\begin{eqnarray}
C_{W,M,m}(u,v,w) &=& \mathbb{P}(\Phi(-W_1) \leq u, 2 \Phi(-m_1) -1 \leq v, 2 \Phi(-M_1) \leq w) \\
&=& \mathbb{P}(1-U \leq u, 1- W \leq v, 1- V \leq w)  \nonumber \\
&=& \mathbb{P}(U \geq 1-u, W \geq 1-v, V \geq 1-w) \nonumber \\
&=& C_{W,M,m}(1,1,1) 
- C_{W,M,m}(1,1,1-v)- C_{W,M,m}(1-u,1,1)- C_{W,M,m}(1,1-w,1) \nonumber\\
&& +\, C_{W,M,m}(1,1-w,1-v) + C_{W,M,m}(1-u,1,1-v) + C_{W,M,m}(1-u,1-w,1) \nonumber\\
&&- \, C_{W,M,m}(1-u,1-w,1-v)\nonumber\\
&=& u+v+w-2 +  C_{W,M,m}(1,1-w,1-v) + C_{W,M,m}(1-u,1,1-v) \nonumber\\
&& + \, C_{W,M,m}(1-u,1-w,1)
- C_{W,M,m}(1-u,1-w,1-v). \nonumber
\end{eqnarray}
%%\end{adjustwidth}

The involutions of the marginal copulas are obtained by substitution in the preceding formulas, respectively, where $u=1$, $v=1$, or $w=1$.

\begin{proof}[Proof of Proposition \ref{P:rho}.]
We recall that the Spearman $\rho$ of the copula $C$ is given by the formula (\cite{Nel} Theorem 5.1.6)%\vspace{-6pt}
\begin{equation}\label{E:rho}
\rho(C) = 12 \int_0^1 \int_0^1 C(u,v) du dv - 3.
\end{equation}

Due to the involution ing{(\ref{E:Involution_2}),}%%%: we revised the format of equation citation, please confirm
 the Spearman $\rho$ of the copula $C_{W,M}$ is the same as the Spearman $\rho$ of the copula $C_{W,m}$.
We observe that
\begin{eqnarray}
\int_0^1 \int_0^1 C_{W,m}(u,w) du dw &=&
\int_0^1 \int_0^{w/2} u du dw + \int_0^1 \int_{w/2}^1 w du dw \\
&& - \, \int_0^1 \int_{w/2}^1 \Phi(2\Phi^{-1}(w/2) - \Phi^{-1}(u)) du dw. \nonumber
\end{eqnarray}

The first two integrals are obvious:
\[ \int_0^1 \int_0^{w/2} u du dw = \int_0^1 \frac{1}{8}w^2 dw = \frac{1}{24}, \]
\[ \int_0^1 \int_{w/2}^1 w du dw = \int_0^1  w\left(1 - \frac{w}{2} \right) dw = \frac{1}{2} - \frac{1}{6}= \frac{1}{3}.\] 

To calculate the third integral, we substitute
\[ u = \Phi(x), \;\; w =2 \Phi(z).\]

We obtain %\vspace{-9pt}
%\begin{adjustwidth}{-\extralength}{0cm}
%\centering %% If there is a figure in wide page, please release command \centering
\begin{eqnarray}
\int_0^1 \int_{w/2}^1 \Phi(2\Phi^{-1}(w/2) - \Phi^{-1}(u)) du dw &=& 2 \int_{-\infty}^0 \int_z^{+\infty} \Phi(2z-x) \phi(x)\phi(z) dx dz \\
&=& 2 \int_{-\infty}^0 \int_z^{+\infty} \int_{-\infty}^{2z-x} \phi(\xi) \phi(x)\phi(z) d\xi dx dz. \nonumber
\end{eqnarray}
%%\end{adjustwidth}

Note that the product $ \phi_3(\xi,x,z)=\phi(\xi) \phi(x)\phi(z)$ is a density of the three-dimensional standard normal distribution $N(0,Id)$, which is a spherical distribution.
Therefore, the above integral can be expressed in terms of the volume of a spherical triangle. We recall that such a volume equals
\begin{equation}\label{E:volume}
 Vol= A+B+C - \pi,
 \end{equation}
 where $A,B,C$ are inner angles  (in radians) of the triangle (see, for example,~\cite{Korn_1968} Formula (1.12.2)  in Section 1.12).
 Note that the cosine of the inner angle of the intersection of two half-spaces is given by  the scalar product of normal vectors:
 \begin{equation}
 \cos(\angle (a_1\xi + b_1 x + c_1 z \leq 0, \, a_2\xi + b_2 x + c_2 z \leq 0) )= - \frac{a_1a_2+b_1b_2+c_1c_2}{\sqrt{a_2^2+b_2^2+c_2^2}\sqrt{a_1^2+b_1^2+c_1^2}}.
 \end{equation}
%\vspace{-23pt}

%\begin{adjustwidth}{-\extralength}{0cm}
%\centering %% If there is a figure in wide page, please release command \centering
\begin{eqnarray}
 \int_{-\infty}^0 \int_z^{+\infty} \int_{-\infty}^{2z-x} \phi(\xi) \phi(x)\phi(z) d\xi dx dz 
 &=& \iiint_{\{z \leq 0,z\leq x,\xi\leq 2z-x\}}\phi_3(\xi,x,z) d\xi dx dz \\
 &=& \frac{1}{\mu_{S^2}(S^2)} \mu_{S^2} (\{(\xi,x,z) \in S^2 : z \leq 0, \,z\leq x ,\, \xi \leq 2z-x \} \nonumber \\
 &=& \frac{1}{4\pi} \bigg( \angle (z \leq 0,z-x\leq 0) + \angle (z\leq 0,\xi+x-2z \leq 0) \nonumber \\
&& + \, \angle(z-x\leq 0,\xi+x-2z \leq 0) - \pi \bigg) \nonumber \\
 &=& \frac{1}{4\pi} \left( \frac{3}{4}\pi + \frac{1}{6}\pi + \arccos\left( \frac{\sqrt{6}}{3}\right) - \pi \right) \nonumber \\
 &=& \frac{1}{4\pi} \arccos\left( \frac{\sqrt{6}}{3}\right) - \frac{1}{48}. \nonumber
\end{eqnarray}
%%\end{adjustwidth}

Finally, we {conclude the following:} %%%: We changed the commas between the digits into decimal dots in equation. Please confirm this revision and check all equations.
%\begin{adjustwidth}{-\extralength}{0cm}
%\centering %% If there is a figure in wide page, please release command \centering
\[ \rho(C_{W,m})=12 \left( \frac{1}{3} + \frac{1}{24} - 2 \left( \frac{1}{4\pi} \arccos\left( \frac{\sqrt{6}}{3}\right) - \frac{1}{48} \right) \right) -3
=2 - \frac{6}{\pi} \arccos\left( \frac{\sqrt{6}}{3}\right) \approx 0.82452.\]
%%\end{adjustwidth}

The case of the copula $C_{M,m}$ can be handled in a similar way. We  substitute
\[ v = 2\Phi(y)-1, \;\; w =2 \Phi(z).\]
%\vspace{-40pt}

%\begin{adjustwidth}{-\extralength}{0cm}
%\centering %% If there is a figure in wide page, please release command \centering
\begin{eqnarray}
I &=&
\int_0^1\int_0^1 C_{M,m}(v,w) dvdw\\
&=&
2\int_{-\infty}^0 \int_0^{+\infty} \bigg( \Psi(z,2y,2(y-z)) - \Psi(2z-y,2y,2(y-z)) \bigg) \varphi(y)\varphi(z) dy dz  \nonumber \\
&=& 2 \sum_{k=0}^\infty \int_{-\infty}^0 \int_0^{+\infty} \bigg( \Phi(z-2k(y-z)) - \Phi(z-2y-2k(y-z))\nonumber \\
&& - \, \Phi(2z-y-2k(y-z)) + \Phi(2z-3y-2k(y-z)) \bigg) \varphi(y)\varphi(z) dy dz \nonumber \\
&=& 2 \sum_{k=0}^\infty \left( 
\int_{-\infty}^0 \int_0^{+\infty} \int_{-\infty}^{(2k+1)z-2ky}\varphi_3(\xi,y,z) d\xi dy dz 
-\int_{-\infty}^0 \int_0^{+\infty} \int_{-\infty}^{(2k+1)z-(2k+2)y}\varphi_3(\xi,y,z) d\xi dy dz \right. \nonumber \\
&& \left. - \, \int_{-\infty}^0 \int_0^{+\infty} \int_{-\infty}^{(2k+2)z-(2k+1)y}\varphi_3(\xi,y,z) d\xi dy dz 
+\int_{-\infty}^0 \int_0^{+\infty} \int_{-\infty}^{(2k+2)z-(2k+3)y}\varphi_3(\xi,y,z) d\xi dy dz 
\right). \nonumber
\end{eqnarray}
%%\end{adjustwidth}

Next, we observe that for $a,b \geq 0$, we have
%\begin{adjustwidth}{-\extralength}{0cm}
%\centering %% If there is a figure in wide page, please release command \centering
\begin{eqnarray}
I_{a,b}&=&
\int_{-\infty}^0 \int_0^{+\infty} \int_{-\infty}^{a z- b y}\varphi_3(\xi,y,z) d\xi dy dz\\
 &=&
\frac{1}{4\pi} \bigg( \angle (z \leq 0, -y \leq 0) + \angle (z\leq 0,\xi -az +by  \leq 0) 
+ \angle (-y\leq 0,\xi -az +by  \leq 0)  - \pi \bigg)
\nonumber \\
&=& \frac{1}{4\pi} \bigg( \frac{1}{2}\pi + \arccos \left( \frac{a}{\sqrt{1+a^2+b^2}}\right)+ \arccos \left( \frac{b}{\sqrt{1+a^2+b^2}}\right) - \pi \bigg) \nonumber \\
&=& \frac{1}{4\pi} \bigg( \arccos \left( \frac{a}{\sqrt{1+a^2+b^2}}\right)+ \arccos \left( \frac{b}{\sqrt{1+a^2+b^2}}\right) - \frac{1}{2}\pi \bigg). \nonumber
\end{eqnarray}
%%\end{adjustwidth}

Note that $I_{a,b}=I_{b,a}$. Therefore, putting respectively $n=2k$ and $n=2k+1$, we obtain
\begin{eqnarray}
I &=& 2 \sum_{k=0}^\infty \bigg( I_{2k+1,2k} - I_{2k+2,2k+1} - I_{2k+2,2k+1} + I_{2k+3,2k+2} \bigg)\\
&=& 2 \sum_{n=0}^\infty (-1)^n \bigg( I_{n+1,n} -I_{n+2,n+1} \bigg) \nonumber \\
&=&  \sum_{n=0}^\infty (-1)^n  \frac{1}{2\pi} \bigg( 
\arccos \left( \frac{n}{\sqrt{2(n^2+n+1)}}\right)
+ \arccos \left( \frac{n+1}{\sqrt{2(n^2+n+1)}}\right) \nonumber \\
&& - \, \arccos \left( \frac{n+1}{\sqrt{2(n^2+3n+3)}}\right)
- \arccos \left( \frac{n+2}{\sqrt{2(n^2+3n+3)}}\right)
\bigg). \nonumber 
\end{eqnarray}

Finally, we conclude the following: %\vspace{-12pt}
%\begin{adjustwidth}{-\extralength}{0cm}
%\centering %% If there is a figure in wide page, please release command \centering
\begin{eqnarray}
\rho(C_{M,m})&=&12 I -3 = \frac{6}{\pi} \sum_{n=0}^\infty (-1)^n \bigg(
\arccos \left( \frac{n}{\sqrt{2(n(n+1)+1)}}\right)
+ \arccos \left( \frac{n+1}{\sqrt{2(n(n+1)+1)}}\right) \nonumber \\
&& - \, \arccos \left( \frac{n+1}{\sqrt{2((n+1)(n+2)+1)}}\right)
- \arccos \left( \frac{n+2}{\sqrt{2((n+1)(n+2)+1)}}\right)
\bigg) \\
&\approx & 0.80649. \nonumber
\end{eqnarray}
%%\end{adjustwidth}
\end{proof}

%\pagebreak

\section{Application of Results: Double-Barrier Option Pricing}\label{S:Application}

In this section, we present one of the possible applications of the derived results. We consider probability space $(\Omega, \mathcal{F}, \mathbb{P})$ and the Black--Scholes--Merton market model~\cite{BS73} or~\cite{Musiela_Rutkowski_2005} in which the asset price $S_t$ dynamic follows the geometric Brownian motion, i.e.,
$$dS_t = \mu S_t dt + \sigma S_t dW^1_t,$$
where $\mu \in \mathbb{R}$, $\sigma > 0$, and $(W^1_t)_{t \ge 0}$ is a Wiener process.
Additionally, the bank account $(B_t)_{t \ge 0}$ is described with the following differential equation:
$$dB_t = rB_t dt,~~ B_0=1.$$

The solution to this system of equations is well known in the literature, and we present only a final result without proof. The reader can check that
$$S_t = S_0e^{\sigma W_t}, \mbox{where } W_t={(\frac{\mu}{\sigma}-\frac{1}{2}\sigma)t+ W^1_t},$$
and
$$B_t = e^{rt}$$
indeed satisfy the above stochastic differential equations.

To valuate in this model a derivative %at time $t<T$ 
with a payoff $X$, which is $\mathcal{F}_T$ measurable ($(\mathcal{F}_t)_{t\ge0}$ is the natural filtration of $W^1_t$), we need to calculate the expected value of the form 
$$V= \mathbb{E}_{\mathbb{P}^*}\bigg(\frac{X}{B_T}\bigg)= e^{-rT} \mathbb{E}_{\mathbb{P}^*}(X),$$
where $\mathbb{P}^*$ is {so-called risk neutral measure,} i.e., the
probability measure equivalent to $\mathbb{P}$ and described by the following Radon--Nikodym derivative
$$\frac{d\mathbb{P}^*}{d\mathbb{P}}=
\exp\left({\frac {r-\mu }{\sigma }}W^1_{T}-{\frac {1}{2}}({\frac {r-\mu }{\sigma }})^{2}T\right).$$

{Note that,} 
 for $t \leq T$, $\,W^2_t=W^1_t-\frac{r - \mu}{\sigma}t$ coincides with a Wiener process in $\mathbb{P}^*$ measure.

For example the Black--Scholes formula for the price of an European call option is
\[C_{BS}(S_0,K,T,r,\sigma) = e^{-rT}\mathbb{E}_{\mathbb{P}^*}((S_T-K)^+).\]

Notice that 
\[ W_t= \left( \frac{r}{\sigma} -\frac{1}{2} \sigma \right) t + W^2_t.\]
%$$S_t = S_0e^{\sigma W_t}=S_0e^{(\mu-\frac{1}{2}\sigma^2)t+\sigma W^1_t}=S_0e^{(\mu-\frac{1}{2}\sigma^2)t+\sigma W^2_t+(r-\mu)t}=
%S_0e^{(r-\frac{1}{2}\sigma^2)t+\sigma W^2_t}$$

This means that {under the risk neutral} measure $\mathbb{P}^*$, the stochastic differential equation of $S_t=S_0e^{\sigma W_t}$ is of the form
\begin{equation}
dS_t = r S_t dt + \sigma S_t dW^2_t.
\end{equation}

We continue to apply the Girsanov theorem. Let $\mathbb{Q} $ be the probability measure given by 
the Radon--Nikodym derivative of the form
\[ \frac{d \mathbb{Q}}{d\mathbb{P}^*}
= \exp \left( -\left( \frac{r}{\sigma} -\frac{1}{2} \sigma \right) W^2_T
- \frac{1}{2} \left( \frac{r}{\sigma} -\frac{1}{2} \sigma \right)^2 T \right).
\]

In the measure $Q$, for $t \leq T$,  $\; W_t = (\ln S_t -\ln S_0)/\sigma$ coincides with a Wiener process. {The reader is referred to~\cite{Oksendal_2003} Theorem 8.6.4 for more details.}
As in Section \ref{S:Main}, we will denote by $M_t$ and $m_t$ the running maxima and minima of $W_t$.

Note that the inverse Radon--Nikodym derivative is given by the formula
\[\frac{d\mathbb{P}^*}{d \mathbb{Q}}
= 
\exp\bigg( \bigg(\frac{r}{\sigma} -\frac{1}{2} \sigma\bigg) W_T
- \frac{1}{2} \bigg(\frac{r}{\sigma} -\frac{1}{2} \sigma\bigg)^2 T \bigg).\]

Now, we can move on to the valuation of the double-barrier options. 

A double-barrier option is a derivative whose payoff is path-dependent. There are few possibilities, like starting between arbitrary values $L < S_0 < U$, and the option becomes worthless whenever $S_t$ crosses at least one of these two barriers.
To be precise, we restrict ourselves to the double-barrier European call option with the payoff $X$ at $T$ of the form
\begin{equation}\label{payoff}
X = (S_T-K)^+\mathds{1}_{L < \min_{0 \le t \le T} S_t} \mathds{1}_{\max_{0 \le t \le T} S_t< U},
\end{equation}
where $\mathds{1}_{\mathcal{A}}$ is the indicator of the event $\mathcal{A}$, and $0 < L < K,S_0 < U$.

First of all, notice that
$$L < \min_{0 \le s \le t}S_s$$
is equivalent to
$$
m_t \geq z, \;\;\; z= \frac{1}{\sigma}\ln \left( \frac{L}{S_0} \right).
$$

Analogously,
$$ \max_{0 \le s \le t} S_s< U $$
is equivalent to
$$
M_t \leq y, \;\;\; y= \frac{1}{\sigma}\ln \left( \frac{U}{S_0} \right).
$$
%\pagebreak

Now, we are ready to move on toward the final price of the derivative.
\begin{Theorem}
The price of a European double-barrier option maturing at time $T$ with a strike $K$ and barriers $L,U$ satisfying
$0 < L < K,S_0 < U$ with the payoff (\ref{payoff}) of the form
\[ C_{2B}(S_0,L,U,K,T,r,\sigma) =\sum_{k=-\infty}^{+\infty} C_{2B*}(k), \]
where %\vspace{-12pt}
%\begin{adjustwidth}{-\extralength}{0cm}
%\centering %% If there is a figure in wide page, please release command \centering
\begin{eqnarray}
C_{2B*}(k) &=&
\exp\bigg( \bigg(\frac{r}{\sigma} -\frac{1}{2} \sigma\bigg) 2k(y-z)  \bigg) 
 %\times
 \left( C_{BS}(e^{2 k \sigma (y-z)} S_0 , K, T, r, \sigma) \right.\\
&& \left. \; -  C_{BS}(e^{2 k \sigma (y-z)} S_0 , U, T, r, \sigma) 
- e^{-rT}(U-K) \Phi\left(\frac{2k(y-z)-y +(r/\sigma-\sigma/2)T}{\sqrt{T}}\right) \right)\nonumber\\
&-&
  \exp\bigg( \bigg(\frac{r}{\sigma} -\frac{1}{2} \sigma\bigg) (2k(y-z)+2y)  \bigg) 
 \left( C_{BS}(e^{2 \sigma (k(y-z) +y) } S_0 , K, T, r, \sigma) \right. \nonumber \\
&& \left. \; -  C_{BS}(e^{2  \sigma (k(y-z)+y)} S_0 , U, T, r, \sigma)  
 - e^{-rT}(U-K) \Phi\left(\frac{2k(y-z)+y+(r/\sigma-\sigma/2)T}{\sqrt{T}}\right) \right) , \nonumber \\
&& y= \frac{1}{\sigma} \ln \left( \frac{U}{S_0} \right), \;\;\;\; z= \frac{1}{\sigma} \ln \left( \frac{L}{S_0} \right). \nonumber
\end{eqnarray}
%%\end{adjustwidth}
\end{Theorem}

%\pagebreak

Note that since, for $k$ tending to $\pm\infty$, $\, C_{2B}(k)$  are converging to 0, for practical purposes, it is enough to consider only a couple of terms with small $|k|$ values. {The similar formulas were provided by Haug (\cite{Haug_1999}) and by Barker (\cite{Barker_2007}).}

\begin{proof}
The price of a European double-barrier option is given by a formula:
\begin{eqnarray}
C_{2B}&=&C_{2B}(S_0,L,U,K,T,r,\sigma)\\
&=& e^{-rT}\mathbb{E}_{P^\ast}\left( (S_T-K)^+\mathds{1}_{L < \min_{0 \le s \le T} S_s} \mathds{1}_{\max_{0 \le s \le T} S_s< U} \right ) \nonumber\\
&=& 
e^{-rT}\mathbb{E}_{Q}\left( \frac{dP^\ast}{dQ} (S_T-K)^+\mathds{1}_{m_T \geq z} \mathds{1}_{M_T \leq y} \right ) \nonumber\\
&=& \mathbb{E}_{Q} \left(  H(W_T) \mathds{1}_{m_T \geq z} \mathds{1}_{M_T \leq y} \right). \nonumber
\end{eqnarray}
where
\[ H(W_T) = e^{-rT} \frac{dP^\ast}{dQ} (S_T-K)^+\mathds{1}_{W_T \leq y},\]
and %% \vspace{-12pt}
%\begin{adjustwidth}{-\extralength}{0cm}
%\centering %% If there is a figure in wide page, please release command \centering
\begin{eqnarray}
H(\xi)&=& e^{-rT} \exp\bigg( \bigg(\frac{r}{\sigma} -\frac{1}{2} \sigma\bigg) \xi
- \frac{1}{2}\bigg( \bigg(\frac{r}{\sigma} -\frac{1}{2} \sigma\bigg)^2 T \bigg)\left(S_0 exp(\sigma \xi) -K\right)^+ \mathds{1}_{\xi \leq y}\\
&=& e^{-rT} \exp\bigg( \bigg(\frac{r}{\sigma} -\frac{1}{2} \sigma\bigg) \xi
- \frac{1}{2}\bigg( \bigg(\frac{r}{\sigma} -\frac{1}{2} \sigma\bigg)^2 T \bigg)
\left( \left(S_0 e^{\sigma \xi} -K\right)^+ - \left(S_0 e^{\sigma \xi} -U\right)^+ -(U-K)\mathds{1}_{\xi \geq y} \right) . \nonumber
\end{eqnarray}
%%\end{adjustwidth}

Next, we apply Theorem \ref{T:F}. Basing on the fact that the density of the joint distribution of $(W_T,M_T,m_T)$, $\;\; F_T(\xi,\eta, \theta)$, is vanishing outside the set
\[ \{(\xi,\eta,\theta) \in \mathbb{R}^3:\;  \theta \leq -\xi^-, \; \xi^+ \leq \eta\}, \;\;\; \xi^\pm = \max(0, \pm \xi) ,\]
we obtain

%\begin{adjustwidth}{-\extralength}{0cm}
%\centering %% If there is a figure in wide page, please release command \centering
\begin{eqnarray}
C_{2B}
&=& \int_{-\infty}^{+\infty} H(\xi) \int_{\xi^+}^y \int_z^{-\xi^-} \frac{\partial^3 F_T}{\partial \xi \partial \eta \partial \theta}(\xi, \eta, \theta) d \theta d \eta d \xi \\
&=& \int_{-\infty}^{+\infty} H(\xi) \bigg( \partial_\xi F_T(\xi,y, -\xi^-) + \partial_\xi F_T(\xi,\xi^+,z) 
 - \partial_\xi F_T(\xi,y,z) \nonumber \\
&&  - \partial_\xi F_T(\xi, \xi^+, -\xi^-) \bigg) d\xi \nonumber \\
&=& \int_{-\infty}^{+\infty} H(\xi)  \bigg( \varphi_T(\xi) - \varphi_T(\xi -2y) \nonumber \\
&&+\sum_{k=0}^\infty \varphi_T(\xi - 2(k+1) (y-z))
-\sum_{k=0}^\infty \varphi_T(\xi - 2(k+2) y+ 2(k+1)z) \nonumber \\
&& - \sum_{k=0}^\infty \varphi_T(-\xi - 2k y+ 2(k+1)z)
+ \sum_{k=0}^\infty \varphi_T(-\xi - 2(k+1)(y-z)) \bigg) d\xi \nonumber \\
&=& \int_{-\infty}^{+\infty} H(\xi) \left(
\sum_{k=0}^\infty \varphi_T(\xi - 2k (y-z))
-\sum_{k=0}^\infty \varphi_T(\xi - 2(k+1) y+ 2kz) \right. \nonumber \\
&&\left. - \sum_{k=0}^\infty \varphi_T(\xi + 2k y - 2(k+1)z)
+ \sum_{k=0}^\infty \varphi_T(\xi + 2(k+1)(y-z))\right) d\xi \nonumber \\
&=&\int_{-\infty}^{+\infty} H(\xi) \left(
\sum_{k=-\infty}^\infty \varphi_T(\xi - 2k (y-z))
-\sum_{k=-\infty}^\infty \varphi_T(\xi - 2(k+1) y+ 2kz) \right) d\xi, \nonumber 
\end{eqnarray}
%%\end{adjustwidth}

where $\varphi_T$ denotes the density of the normal distribution $N(0,T)$, and
\[\varphi_T(x)= \frac{1}{\sqrt{2T \pi}} \exp\left( - \frac{x^2}{2T} \right).\] 

{Note that the formula in the last brackets may be obtained as well from
the L{\'e}vy Triple Law; see Corollary \ref{C:LTL} or~\cite{Schilling_2014} Theorem 6.18.}

Next, we return to the probability $\mathbb{P}^\ast$.

% \vspace{-12pt}
%\begin{adjustwidth}{-\extralength}{0cm}
%\centering %% If there is a figure in wide page, please release command \centering
\begin{eqnarray}
C_{2B}
&=&
\sum_{k=-\infty}^\infty  \int_{-\infty}^{+\infty} H(\xi)  \varphi_T(\xi - 2k (y-z)) d\xi
-\sum_{k=-\infty}^\infty \int_{-\infty}^{+\infty} H(\xi)  \varphi_T(\xi - 2(k+1) y+ 2kz))  d\xi \nonumber \\
&=&
\sum_{k=-\infty}^\infty  \int_{-\infty}^{+\infty} H(\xi + 2k (y-z) ) \varphi_T(\xi ) d\xi
-\sum_{k=-\infty}^\infty \int_{-\infty}^{+\infty} H(\xi + 2(k+1) y- 2kz)  \varphi_T(\xi))  d\xi \nonumber \\
&=&
\sum_{k=-\infty}^\infty  \mathbb{E}_Q \bigg( H(W_T+ 2k (y-z) )  \bigg)
-\sum_{k=-\infty}^\infty \mathbb{E}_Q \bigg( H(W_T+ 2(k+1) y- 2kz) \bigg) 
\end{eqnarray}
%%\end{adjustwidth}

Since
% \vspace{-12pt}
%\begin{adjustwidth}{-\extralength}{0cm}
%\centering %% If there is a figure in wide page, please release command \centering
\begin{eqnarray}
H(W_T + \Delta) &=&
e^{-rT} \exp\bigg( \bigg(\frac{r}{\sigma} -\frac{1}{2} \sigma\bigg) (W_T+\Delta )
- \frac{1}{2}\bigg( \bigg(\frac{r}{\sigma} -\frac{1}{2} \sigma\bigg)^2 T \bigg) \\
&& \times
\left( \left(S_0 e^{\sigma (W_T+\Delta)} -K\right)^+ - \left(S_0 e^{\sigma (W_T+\Delta)} -U\right)^+ -(U-K)\mathds{1}_{W_T \geq y- \Delta} \right) \nonumber \\
&=& \frac{d \mathbb{P}^\ast}{d \mathbb{Q}} 
e^{-rT} \exp\bigg( \bigg(\frac{r}{\sigma} -\frac{1}{2} \sigma\bigg) \Delta 
 \bigg) \nonumber \\
 && \times
\left( \left(S_0 e^{\sigma \Delta} e^{\sigma W_T} -K\right)^+ - \left(S_0 e^{\sigma \Delta} e^{\sigma W_T} -U\right)^+ -(U-K)\mathds{1}_{W_T \geq y- \Delta} \right), \nonumber
\end{eqnarray}
%%\end{adjustwidth}
we have

%\begin{adjustwidth}{-\extralength}{0cm}
%\centering %% If there is a figure in wide page, please release command \centering
\begin{eqnarray}
\mathbb{E}_Q(H(W_T + \Delta) )&=&
\mathbb{E}_{P^\ast}\left( e^{-rT} \exp\bigg( \bigg(\frac{r}{\sigma} -\frac{1}{2} \sigma\bigg) \Delta 
 \bigg) \right.\\
&& \times \left. \left( \left(S_0 e^{\sigma \Delta} e^{\sigma W_T} -K\right)^+ 
- \left(S_0 e^{\sigma \Delta} e^{\sigma W_T} -U\right)^+ 
-(U-K)\mathds{1}_{W_T \geq y- \Delta} \right) \right) \nonumber\\
&=&
\exp\bigg( \bigg(\frac{r}{\sigma} -\frac{1}{2} \sigma\bigg) \Delta  \bigg) \nonumber \\
&& \times
 \left( C_{BS}(e^{\sigma \Delta} S_0 , K, T, r, \sigma) -  C_{BS}(e^{\sigma \Delta} S_0 , U, T, r, \sigma) \right. \nonumber \\
&& \; \left. - e^{-rT}(U-K) \Phi\left(\frac{\Delta-y +(r/\sigma-\sigma/2)T}{\sqrt{T}}\right) \right)\nonumber
\end{eqnarray}
%%\end{adjustwidth}

{To conclude, for the proof, it is enough to put, respectively, $\Delta=2k(y-z)$ or $\Delta = 2k(y-z)+2y$}.
\end{proof}

{The alternative method of pricing the double-barrier options is via simulations in the risk-neutral measure. In the {Table} %%%: we revised the format as table citation, please confirm
 \ref{Tabela}, we present the results. The simulation is provided by an approximation of the Wiener process by a Random Walk process with {100,000} 
 steps. We iterated the procedure 1534 times. The analytic formula for the price is approximated by the sum of the seven terms $C_{2B*}(k)$, where $k=-3, \dots , 3$.}

\begin{table}[ht]
\caption{{The} %%%: 1. we revised the format as Table 1, and please check the citation in main text as well 2. We changed the commas between the digits into decimal dots. Please confirm this revision. 3. we removed vertical line, please confirm and check
 comparison of simulated (S) and approximated (A) values of the two-barrier call options with varying strike $K$ and risk-free interest rate $r$ (in \%).
The initial value of the stock, barriers, volatility, and time to maturity are fixed: $S_0=100$, $L=80$, $U=130$, $\sigma=0{.}5$, and $T=1/4$.   }
\label{Tabela}
\begin{tabular*}{\textwidth}{ccccccccccccc}
\hline
\boldmath{$r$} & \multicolumn{2}{c}{\textbf{{1} 
}} & \multicolumn{2}{c}{\textbf{3}} & \multicolumn{2}{c}{\textbf{6}} & \multicolumn{2}{c}{\textbf{9}}  & \multicolumn{2}{c}{\textbf{12.5}} \\
\hline
\boldmath{$K$}	&	\textbf{S}	&	\textbf{A}	&	\textbf{S}	&	\textbf{A}	&	\textbf{S}	&	\textbf{A}	&	\textbf{S}	&	\textbf{A}	&	\textbf{S}	&	\textbf{A}	\\
\hline 	
80	&	7.318	&	7.373	&	7.434	&	7.39	&	7.486	&	7.412	&	7.556	&	7.428	&	7.507	&	7.439	\\
81	&	6.981	&	7.037	&	7.095	&	7.055	&	7.149	&	7.078	&	7.219	&	7.095	&	7.178	&	7.109	\\
82	&	6.645	&	6.703	&	6.757	&	6.722	&	6.813	&	6.745	&	6.883	&	6.764	&	6.849	&	6.779	\\
83	&	6.31	&	6.371	&	6.42	&	6.39	&	6.477	&	6.415	&	6.547	&	6.435	&	6.521	&	6.452	\\
84	&	5.977	&	6.042	&	6.086	&	6.062	&	6.144	&	6.088	&	6.214	&	6.109	&	6.195	&	6.127	\\
85	&	5.649	&	5.718	&	5.755	&	5.738	&	5.814	&	5.765	&	5.884	&	5.787	&	5.871	&	5.807	\\
86	&	5.327	&	5.4	&	5.432	&	5.42	&	5.49	&	5.447	&	5.557	&	5.47	&	5.55	&	5.491	\\
87	&	5.012	&	5.087	&	5.113	&	5.108	&	5.173	&	5.135	&	5.239	&	5.159	&	5.233	&	5.182	\\
88	&	4.703	&	4.782	&	4.803	&	4.802	&	4.862	&	4.83	&	4.926	&	4.855	&	4.925	&	4.878	\\
89	&	4.403	&	4.484	&	4.501	&	4.504	&	4.56	&	4.533	&	4.62	&	4.557	&	4.622	&	4.582	\\
90	&	4.113	&	4.194	&	4.207	&	4.215	&	4.265	&	4.243	&	4.323	&	4.268	&	4.326	&	4.293	\\
91	&	3.834	&	3.913	&	3.924	&	3.934	&	3.98	&	3.962	&	4.035	&	3.987	&	4.04	&	4.012	\\
92	&	3.568	&	3.642	&	3.654	&	3.662	&	3.707	&	3.69	&	3.757	&	3.715	&	3.764	&	3.741	\\
93	&	3.311	&	3.38	&	3.395	&	3.4	&	3.446	&	3.428	&	3.492	&	3.453	&	3.497	&	3.478	\\
94	&	3.065	&	3.129	&	3.146	&	3.148	&	3.194	&	3.175	&	3.24	&	3.2	&	3.244	&	3.225	\\
95	&	2.828	&	2.888	&	2.908	&	2.907	&	2.955	&	2.933	&	2.998	&	2.957	&	3.002	&	2.982	\\
96	&	2.601	&	2.658	&	2.678	&	2.676	&	2.725	&	2.701	&	2.768	&	2.725	&	2.77	&	2.749	\\
97	&	2.388	&	2.438	&	2.461	&	2.455	&	2.504	&	2.48	&	2.546	&	2.503	&	2.549	&	2.527	\\
98	&	2.183	&	2.229	&	2.255	&	2.246	&	2.296	&	2.269	&	2.334	&	2.291	&	2.337	&	2.315	\\
99	&	1.988	&	2.031	&	2.056	&	2.047	&	2.097	&	2.07	&	2.136	&	2.091	&	2.137	&	2.113	\\
100	&	1.806	&	1.844	&	1.871	&	1.859	&	1.906	&	1.881	&	1.945	&	1.901	&	1.949	&	1.923	\\
101	&	1.633	&	1.668	&	1.696	&	1.682	&	1.73	&	1.703	&	1.763	&	1.722	&	1.768	&	1.742	\\
102	&	1.47	&	1.503	&	1.529	&	1.516	&	1.561	&	1.535	&	1.595	&	1.553	&	1.598	&	1.573	\\
103	&	1.318	&	1.348	&	1.374	&	1.361	&	1.403	&	1.378	&	1.434	&	1.395	&	1.441	&	1.414	\\
104	&	1.176	&	1.204	&	1.229	&	1.215	&	1.256	&	1.232	&	1.284	&	1.248	&	1.291	&	1.265	\\
105	&	1.045	&	1.07	&	1.095	&	1.08	&	1.119	&	1.096	&	1.146	&	1.11	&	1.151	&	1.126	\\
106	&	0.923	&	0.946	&	0.972	&	0.956	&	0.993	&	0.97	&	1.017	&	0.983	&	1.023	&	0.998	\\
107	&	0.813	&	0.832	&	0.857	&	0.84	&	0.877	&	0.853	&	0.9	&	0.865	&	0.904	&	0.879	\\
\hline
%\end{tabularx}}
\end{tabular*}
\end{table}

\begin{table}[ht]%\ContinuedFloat
\small
\caption{{\em Cont.}}
\begin{tabular*}{\textwidth}{ccccccccccccc}
\hline
\boldmath{$r$} & \multicolumn{2}{c}{\textbf{{1} 
}} & \multicolumn{2}{c}{\textbf{3}} & \multicolumn{2}{c}{\textbf{6}} & \multicolumn{2}{c}{\textbf{9}}  & \multicolumn{2}{c}{\textbf{12.5}} \\
\hline
\boldmath{$K$}	&	\textbf{S}	&	\textbf{A}	&	\textbf{S}	&	\textbf{A}	&	\textbf{S}	&	\textbf{A}	&	\textbf{S}	&	\textbf{A}	&	\textbf{S}	&	\textbf{A}	\\
\hline		
108	&	0.712	&	0.727	&	0.753	&	0.735	&	0.77	&	0.746	&	0.792	&	0.757	&	0.795	&	0.769	\\
109	&	0.621	&	0.631	&	0.658	&	0.638	&	0.673	&	0.648	&	0.693	&	0.658	&	0.696	&	0.669	\\
110	&	0.54	&	0.544	&	0.573	&	0.55	&	0.585	&	0.559	&	0.604	&	0.568	&	0.606	&	0.578	\\
111	&	0.465	&	0.465	&	0.496	&	0.471	&	0.507	&	0.479	&	0.524	&	0.487	&	0.525	&	0.495	\\
112	&	0.395	&	0.394	&	0.424	&	0.399	&	0.436	&	0.406	&	0.452	&	0.413	&	0.452	&	0.421	\\
113	&	0.338	&	0.331	&	0.362	&	0.335	&	0.37	&	0.341	&	0.387	&	0.347	&	0.387	&	0.354	\\
114	&	0.285	&	0.275	&	0.308	&	0.279	&	0.315	&	0.284	&	0.325	&	0.289	&	0.328	&	0.295	\\
115	&	0.237	&	0.226	&	0.257	&	0.229	&	0.266	&	0.233	&	0.275	&	0.237	&	0.273	&	0.242	\\
116	&	0.195	&	0.183	&	0.213	&	0.185	&	0.22	&	0.189	&	0.23	&	0.193	&	0.228	&	0.197	\\
117	&	0.159	&	0.146	&	0.175	&	0.148	&	0.18	&	0.151	&	0.188	&	0.154	&	0.189	&	0.157	\\
118	&	0.127	&	0.114	&	0.141	&	0.116	&	0.146	&	0.118	&	0.152	&	0.12	&	0.151	&	0.123	\\
119	&	0.097	&	0.087	&	0.111	&	0.089	&	0.115	&	0.091	&	0.121	&	0.092	&	0.119	&	0.095	\\
120	&	0.073	&	0.065	&	0.085	&	0.066	&	0.088	&	0.068	&	0.093	&	0.069	&	0.093	&	0.071	\\
121	&	0.054	&	0.047	&	0.063	&	0.048	&	0.065	&	0.049	&	0.069	&	0.05	&	0.07	&	0.051	\\
122	&	0.038	&	0.033	&	0.046	&	0.034	&	0.047	&	0.034	&	0.049	&	0.035	&	0.049	&	0.036	\\
123	&	0.025	&	0.022	&	0.032	&	0.022	&	0.033	&	0.023	&	0.034	&	0.023	&	0.033	&	0.024	\\
124	&	0.014	&	0.014	&	0.02	&	0.014	&	0.021	&	0.014	&	0.023	&	0.015	&	0.021	&	0.015	\\
125	&	0.007	&	0.008	&	0.011	&	0.008	&	0.012	&	0.008	&	0.014	&	0.008	&	0.013	&	0.009	\\
126	&	0.002	&	0.004	&	0.005	&	0.004	&	0.005	&	0.004	&	0.008	&	0.004	&	0.007	&	0.004	\\
127	&	0	&	0.002	&	0.001	&	0.002	&	0.002	&	0.002	&	0.008	&	0.002	&	0.004	&	0.002	\\
128	&	0	&	0	&	0	&	0.001	&	0	&	0.001	&	0.001	&	0.001	&	0.001	&	0.001	\\
129	&	0	&	0	&	0	&	0	&	0	&	0	&	0	&	0	&	0	&	0	\\
130	&	0	&	0	&	0	&	0	&	0	&	0	&	0	&	0	&	0	&	0	\\
\hline
\end{tabular*}
\end{table}

%\pagebreak 

%\begin{adjustwidth}{-\extralength}{0cm}
%\printendnotes[custom] % Un-comment to print a list of endnotes

%\reftitle{References}

% Please provide either the correct journal abbreviation (e.g. according to the “List of Title Word Abbreviations” http://www.issn.org/services/online-services/access-to-the-ltwa/) or the full name of the journal.
% Citations and References in Supplementary files are permitted provided that they also appear in the reference list here. 

%=====================================
% References, variant A: external bibliography
%=====================================
%\bibliography{your_external_BibTeX_file}

%=====================================
% References, variant B: internal bibliography
%=====================================

\end{document}